\documentclass{amsart}
\title[
LLN for finite-range dependent 
random matrices]{A law of large numbers 
for finite-range dependent random matrices} 
\author{Greg Anderson}
\address{School of Mathematics, University of Minnesota, 
206 Church St. SE,
Minneapolis, MN 55455}
\author{Ofer Zeitouni}  
\thanks{ O.Z. was partially supported by NSF grant number
DMS-0503775.}
\address{School  of Mathematics, University of Minnesota, 
206 Church St. SE,
Minneapolis, MN 55455, and
Departments of Mathematics and of EE, Technion, Haifa 32000, Israel}
\date{September 10, 2006. Revised October 2, 2007}

\newcommand{\norm}[1]{{\Vert#1\Vert}}

\newtheorem{Proposition}[subsection]{Proposition}
\newtheorem{Theorem}[subsection]{Theorem}
\newtheorem{Assumption}[subsubsection]{Assumption}
\newtheorem{Lemma}[subsubsection]{Lemma}

\newcommand{\WWW}{{\mathcal{W}}}
\newcommand{\one}{{\mathbf{1}}}
\newcommand{\CC}{{\mathbb{C}}}
\newcommand{\ibold}{{\mathbf{i}}}
\newcommand{\jbold}{{\mathbf{j}}}
\newcommand{\NN}{{\mathbb{N}}}
\newcommand{\ZZ}{{\mathbb{Z}}}
\newcommand{\RR}{{\mathbb{R}}}
\newcommand{\E}{{\mathbb{E}}}
\DeclareMathOperator{\trace}{{\mathrm{trace}}}
\newcommand{\III}{{\mathcal{I}}}
\begin{document}

\begin{abstract}
We consider 
random hermitian 
matrices in which distant above-diagonal entries
are independent but nearby entries may be correlated. 
We find the limit of the empirical
distribution of eigenvalues by combinatorial methods.
We also prove that the limit has algebraic Stieltjes
transform by an argument
based on dimension theory of noetherian local rings.
\end{abstract}
\maketitle

\section{Introduction}
Study of the  empirical distribution of eigenvalues of  random
hermitian (or real symmetric)
matrices has a long history, starting with the seminal work of
Wigner \cite{wigner} and Wishart \cite{wishart}. Except in cases
where the joint distribution of eigenvalues is explicitly known, most
results available are asymptotic in nature and based 
on one of the following approaches: (i) the {\it moment method}, i.~e.,
evaluation of expectations of traces of powers of the matrix; 
(ii) appropriate recursions 
for the resolvant, as introduced in \cite{pastur-marchenko};
or (iii) the free probability method
(especially the notion of asymptotic freeness) 
originating with  Voiculescu \cite{voi91}.
A good review of the first two approaches can 
be found in \cite{bai}.
For the third, see \cite{voi01}, and for a somewhat 
more combinatorial 
perspective, \cite{speicher}.
These approaches have  been extended to situations in which 
the matrix analyzed neither possesses i.i.d.\ entries above the diagonal
(as in the Wigner case) nor is it the product of matrices with 
i.i.d.\ entries (as in the Wishart case). We mention in particular
the papers \cite{MPK}, \cite{khorunzhy1}, \cite{Sh96} and 
\cite{guionnet-clt} for results pertaining to the model of 
``random band matrices'', all with independent above-diagonal entries.

In our recent work \cite{AZ} we 
studied convergence of the empirical distribution of eigenvalues of
random band matrices, and developed a combinatorial approach,
based on the moment method, to identify the limit (and also to provide
central limit theorems for linear statistics). 
Here we develop the method  further to handle a class of matrices
with local dependence among entries (we postpone the precise definition
of the class to Section \ref{sec-2}).
 To each random matrix of the class
we associate a random band matrix 
with the same limit of empirical distribution of eigenvalues by a 
process of ``Fourier transformation'', thus 
making it possible to describe the limit in terms of 
our previous work (see Theorem~\ref{Theorem:MainResult}).
 We also prove that the Stieltjes transform of the 
limit is algebraic (see Theorem~\ref{Theorem:MainResultBis}),
and  deduce in turn from well-known facts about algebraic functions
a strong regularity result for $\mu$
(see  Theorem~\ref{Theorem:MainResultTer}).
Algebraicity we prove by a 
general ``soft'' (i.~e., nonconstructive) method based on the 
theory of noetherian 
local rings  (see Theorem~\ref{Theorem:BasicFunctionalEquation})  which ought to be applicable to many more random matrix problems.

To get the flavor of our
results, the reader should imagine a Wigner matrix (i.e., an 
$N$-by-$N$ real symmetric random matrix
with i.i.d.\  above-diagonal entries, each of mean $0$ 
and variance $1/N$), on which a local ``filtering'' operation 
is performed: each entry not near the diagonal 
or an edge is replaced by half the sum of its
four neighbors to northeast, southeast, southwest and northwest. At the end of Section \ref{section:Method} (Theorem~\ref{Theorem:MainResult} taken for granted)
we analyze the ``(NE+SE+SW+NW)-filtered Wigner matrix'' described above.
We find that the limit measure
is the free multiplicative
convolution of the semicircle law (density $\propto\one_{|x|<2}\sqrt{4-x^2}$)
and the arcsine law (density $\propto\one_{0<x<2}/\sqrt{x(2-x)}$). The 
appearance in this example of a free multiplicative convolution has a 
simple explanation (see Proposition~\ref{Proposition:NonFluke}).
We also write down the quartic equation satisfied by the Stieltjes 
transform of the limit measure.

Recently
other authors have considered the empirical 
distribution of eigenvalues for matrices with dependent entries,
see \cite{GT},\cite{Cha},\cite{Schenker}. 
Their class of models does not overlap significantly
with ours. In particular, in all these works and unlike
in our model, the limit of the empirical
measure is always the same as that of a semicircle law 
multiplied by a random or deterministic constant. 

Closest to our work is the recent paper by \cite{HLN}, that
builds upon earlier work by \cite{BDM} and \cite{Gir}. They consider Gram
matrices of the form $X_NX_N^*$ where $X_N$ is a sequence of (non-symmetric)
Gaussian matrices obtained by applying a filtering operation to a matrix with
(complex) Gaussian independent entries. They also consider the case
$(X_N+A_N)(X_N+A_N)^*$ with $A_N$ deterministic and Toeplitz. The Gaussian 
assumption allows them to directly approximate the matrix $X_N$ by a unitary
transformation of a Gaussian matrix with independent (but not identically
distributed) entries, to which the results of \cite{Gir} and \cite{AZ} 
apply. 
An advantage of their approach is that they do not need 
to assume finiteness of the
filter; they need only a summability condition.
We note also
that
 the approach in \cite{HLN} and
\cite{BDM} is based on the study of resolvants rather than moments.

We mention now motivation from electrical engineering.
The analysis of the limiting 
empirical distribution of eigenvalues of random matrices
has recently played an important role in the analysis of communication
systems, see \cite{verdu} for an extensive review. 
In particular, when studying
multi-antenna systems, one often makes the (unrealistic) assumption 
that gains between different pairs of antennas are 
uncorrelated. The models 
studied in this work would allow correlation
 between neighboring antenna pairs.
We do not develop this application further here.

The structure of the paper  is as follows. In the next section we
describe the class of matrices we treat, and state our main results,
namely  Theorem~\ref{Theorem:MainResult} (asserting a law of large numbers) and Theorem~\ref{Theorem:MainResultBis} (asserting algebraicity of a Stieltjes transform).
We also prove  Theorem~\ref{Theorem:MainResultTer} (which is essentially folkloric and explains the regularity implied by algebraicity).
In Section \ref{section:Method} we discuss the limit measure in detail,
and in particular write down algebro-integral equations for its Stieltjes transform,
which we call {\em color equations}.
Section 
\ref{section:MainLimit} provides a computation of limits of traces
of powers of the matrices under consideration.
In 
Section \ref{section:BigCompletion} we complete the proof of
Theorem~\ref{Theorem:MainResult}  by a variance computation.
In Section  \ref{section:Soft} we set up the algebraic machinery needed to prove  Theorem \ref{Theorem:MainResultBis}.
We finish proving the theorem in Section \ref{section:SoftApplied} by analyzing the color equations.

\section{Formulation of the main result}
\label{sec-2}
After defining a class of kernels in \S \ref{sec-2.1},
we define in \S  \ref{sec-2.1a} the class of matrices 
dealt with in this paper and in \S \ref{sec-motivation}
describe the subclass of {\em filtered Wigner matrices}.
To each kernel we associate a measure in \S\ref{subsection:MuDef}.
Finally, we state our main results, Theorems~\ref{Theorem:MainResult}
and \ref{Theorem:MainResultBis}.

\subsection{Kernels}
\label{sec-2.1}

\subsubsection{Color space (motivation)}
As in \cite{AZ}, the spatial change in the variance structure
of the entries of a random matrix is captured by an auxiliary
variable, which we call ``color''. The difference between the
derivation in \cite{AZ} and the present paper is that 
we append to color space an additional variable (related
to the local averaging mentioned in the introduction), that
one should think of as a Fourier variable. The precise definition
follows.

\subsubsection{Color space (formal definition)}
Let $C=[0,1]\times S^1$, where $S^1$ is the unit circle in 
the complex plane.
We call $C$ {\em color space}.
We declare $C$ to be a probability space under the product of
 uniform probability measures on $[0,1]$ and $S^1$, denoted $P$.  
We say that a $C$-valued random variable is 
{\em uniformly distributed} if its law is $P$.
In the sequel, all $L^p(C)$ (resp., $L^p(C\times C)$) spaces, $p\geq 1$, are
taken with respect to the measure $P$ (resp.,  $P\times P$).


\subsubsection{The kernel $s$}
We fix a kernel 
$$s:C\times C\rightarrow\RR$$
 which will govern the local covariance structure of our random matrix model. We impose on $s$ the following conditions.
\begin{Assumption}
\label{Assumption:Kernel}$\;$

\begin{enumerate}
\item
$s$ is a nonnegative symmetric function, i.e.
$$s(c,c')=s(c',c)\geq 0.$$
\item $s$ has a Fourier expansion
$$s(c,c')=\sum_{i,j\in \ZZ}s_{ij}(x,y)\xi^i\eta^j\;\;(c=(x,\xi),\;c'=(y,\eta))$$
where all but finitely many of the coefficients
$$s_{ij}:[0,1]\times [0,1]\rightarrow\CC$$
vanish identically. \\
\item There is a finite partition $\III$ of $[0,1]$
into subintervals of positive length such that every coefficient function $s_{ij}$ is constant on every set of the form $I\times J$ with $I,J\in \III$. \\
\item $s$ is nondegenerate: $\norm{s}_{L^1(C\times C)}>0$.
 \end{enumerate}
\end{Assumption}

\subsubsection{The sets $Q^{(N)}_K$}
For each $N\in \NN$ (here and below $\NN$ denotes the set of positive integers) and $K>0$, we define
$$Q^{(N)}_K=
\left\{i\in \{1,\dots,N\}\left| \min_{x\in \partial \III}|i-Nx|>K\right.\right\}$$
where $\partial \III\subset [0,1]$ is the finite set consisting of all endpoints
of all intervals belonging to the family $\III$.

\subsubsection{Remarks}
 We have
\begin{equation}\label{equation:Symmetries}
\bar{s}_{ij}=s_{-i,-j},\;\;\;s_{ji}(y,x)=s_{ij}(x,y)
\end{equation}
because $s$ is real-valued and symmetric.
Assumptions (I,II,III) imply that 
\begin{equation}\label{equation:NoExceptions}
0\leq s\leq \norm{s}_{L^\infty(C\times C)}<\infty
\end{equation}
 holds everywhere (not just $P\times P$-a.e.).

\subsection{The model}
\label{sec-2.1a}
For each $N\in \NN$, let $$X^{(N)}=[X^{(N)}_{ij}]_{i,j=1}^N$$ be  an $N$-by-$N$  random
hermitian matrix. We impose the following conditions,
where $s$ satisfies Assumption~\ref{Assumption:Kernel}.

\begin{Assumption}
\label{Assumption:Model}$\;$

\begin{enumerate}
\item 
\begin{enumerate}
\item $\forall N\in \NN\;\;\; \forall i,j\in\{1,\dots,N\}\;\;\;EX^{(N)}_{ij}=0$.
\item $\displaystyle\forall k\in\NN\;\;\;\sup_{N=1}^\infty \max_{i,j=1}^N E|X_{ij}^{(N)}|^k<\infty$.\\
\end{enumerate}
\item There exists 
$K>0$ such that for all $N\in \NN$,
the following hold:\\
\begin{enumerate}
\item
$\forall i,j\in \ZZ\;\;\;\max(|i|,|j|)>K\Rightarrow s_{ij}\equiv 0$.\\
\item For all nonempty subsets 
$$A,B\subset \{(i,j)\in \{1,\dots,N\}^2\mid 1\leq i\leq j\leq N\}$$
such that 
$$\min_{(i,j)\in A}\min_{(k,\ell)\in B}\max(|i-k|,|j-\ell|)>K,$$
the $\sigma$-fields 
$$\sigma(\{X_{ij}^{(N)}\mid (i,j)\in A\}),\;\;\;\sigma(\{X_{k\ell}^{(N)}\mid (k,\ell)\in B\})$$ are independent.\\
\item $\forall i,j,k,\ell\in Q^{(N)}_K\;\mathrm{s.t.}\;
\min(j-i,\ell-k)>K$\\
$EX_{ij}^{(N)}\overline{X}_{k\ell}^{(N)}=s_{i-k,\ell-j}(i/N,j/N)$.
\end{enumerate}
\end{enumerate}
\end{Assumption}

\subsubsection{The empirical distribution of eigenvalues}
\label{subsubsection:EmpiricalDistribution}
Let 
$$\lambda_1^{(N)}\leq \lambda_2^{(N)}\leq \ldots
\leq \lambda_N^{(N)}$$
denote the eigenvalues of the
hermitian matrix
$X^{(N)}/\sqrt{N}$,
and let
$$L^{(N)}=N^{-1}\sum_{i=1}^N \delta_{\lambda_i^{(N)}}$$
denote 
the corresponding empirical distribution of the eigenvalues.
We are concerned with the limiting behavior of $L^{(N)}$ as $N\rightarrow\infty$.



\subsubsection{Remarks}
 (i) The existence of $K$ satisfying Assumption~\ref{Assumption:Model}(IIa)
is assured by Assumption~\ref{Assumption:Kernel}(IIb).
(ii) Assumption~\ref{Assumption:Model}(IIb) says that 
the on-or-above diagonal entries of $X^{(N)}$
form a finite-range dependent random field, with $K$ a 
bound for the range of
dependence. This explains the title of the paper. (iii) Assumption~\ref{Assumption:Model}(IIc) fixes variances at each site 
and also short-range local correlations for sites in
 sufficiently general position. (iv) None of our assumptions 
rule out the possibility that all matrices $X^{(N)}$
are real. In other words, we can handle hermitian and real 
symmetric cases uniformly under Assumption~\ref{Assumption:Model}. (v)  
The relations (\ref{equation:Color}), which are equivalent to the reality and symmetry of the kernel $s$,
play a key role in our analysis of the limiting behavior of $L^{(N)}$. New methods would be needed were the reality and symmetry conditions to be relaxed.

\subsubsection{Kernels of pure spatial type}
Let $s:[0,1]\times [0,1]\rightarrow \RR$ be a kernel.
If the kernel $\tilde{s}:C\times C\rightarrow \RR$ defined by the formula
$\tilde{s}((x,\xi),(y,\eta))=s(x,y)$
satisfies Assumption~\ref{Assumption:Kernel},
then by abuse of terminology
we say that $s$ is a kernel of {\em pure spatial type} satisfying Assumption~\ref{Assumption:Kernel},
and with the evident modification
of Assumption~\ref{Assumption:Model}(IIc) we can use $s$ to govern
the covariance structure of our model. In the special case
of kernels of  pure spatial type our model essentially contains  the model of \cite{AZ} in the special case in which color space is a finite set.

\subsubsection{Kernels of pure Fourier type}
\label{subsubsection:PureFourierType}
Let $s:S^1\times S^1\rightarrow \RR$ be a kernel.
If the kernel $\tilde{s}:C\times C\rightarrow \RR$ defined by the formula $\tilde{s}((x,\xi),(y,\eta))=s(\xi,\eta)$
satisfies Assumption~\ref{Assumption:Kernel}
with $\III=\{[0,1]\}$, then by abuse of terminology
we say that $s$ is a kernel of {\em pure Fourier type} satisfying Assumption~\ref{Assumption:Kernel}, and with the evident modification
of Assumption~\ref{Assumption:Model}(IIc) we can use $s$ to govern
the covariance structure of our model.
We suggest that the reader focus on the pure Fourier case when first approaching this paper because little would be lost in terms of grasping
the main ideas. The main 
reason for us to work at
a higher level of generality is to make sure
that our theory handles not only ``filtered Wigner matrices'' (for which $| \III|=1$) but also ``filtered Wishart matrices"
(for which $| \III|=2$).
(Here and below $|S|$ denotes the cardinality of a set $S$.)
We then might as well allow $|\III|>2$ as a possibility because there is no gain in simplicity by excluding it.

 \subsection{Filtered Wigner matrices}\label{sec-motivation}
 We describe now a natural class of random matrices 
fitting into the framework of Assumptions~\ref{Assumption:Kernel} and \ref{Assumption:Model}.
This class should be considered  the main motivation for the paper.
Members of the class arise by ``filtering'' Wigner matrices.
The corresponding kernels are of pure Fourier type
and depend in a simple way on the ``filter''.

 \subsubsection{Wigner matrices}
Let
$$\{Y_{ij}\}_{-\infty<i<j<\infty}$$ 
be an i.i.d.\  family of real random variables.
Assume that $Y_{01}$ has absolute moments of all orders.
Assume that $EY_{01}=0$ and $EY_{01}^2=1$.
Put 
$$Y_{ii}=0\;\mbox{for $-\infty<i<\infty$},\;\;
Y_{ij}=Y_{ji}\;\mbox{for $-\infty<j<i<\infty$}.
$$
Then
\begin{equation}\label{equation:YCorrelation}
Y_{ij}=Y_{ji},\;\;\;\;E Y_{ij}=0,\;\;\;E Y_{ij}Y_{k\ell}=
(\delta_{ik}\delta_{j\ell}+\delta_{i\ell}\delta_{jk})(1-\delta_{ij})(1-\delta_{k\ell})
\end{equation}
for all $i,j,k,\ell\in \ZZ$.  
Let $Y^{(N)}=[Y^{(N)}_{ij}]_{i,j=1}^N$ be the $N$-by-$N$ matrix with entries $Y_{ij}$.
Then $Y^{(N)}/\sqrt{N}$  (in the terminology of \cite{AZ}) is a {\em Wigner matrix},
and hence  the empirical distribution of its eigenvalues for $N\rightarrow\infty$ tends  to the semicircle law.
In particular, if the $Y_{ij}$ are standard normal random variables
and one makes a suitable adjustment to the diagonal
of $Y^{(N)}/\sqrt{N}$, the result is a Wigner matrix in the 
standard sense, i.e., a  member of the Gaussian orthogonal ensemble.

\subsubsection{The filter, its Fourier transform, and associated kernel}
Let a {\it filter}
$$h:\ZZ\times\ZZ\to\RR$$
be given, with $H$ denoting its Fourier transform, that is
$$H(\xi,\eta)=\sum_{i,j\in \ZZ} h(i,j)\xi^i\eta^j\;\;\;(\xi,\eta\in S^1).$$
We assume that $h$ does not vanish identically.
We assume that there exists $K>0$ such that
\begin{equation}\label{equation:LittlehSupport}
\max(|i|,|j|)>K/2\Rightarrow h(i,j)=0,
\end{equation}
and hence $H$ is well-defined.
We assume that $h$ satisfies the symmetry condition
\begin{equation}\label{equation:LittlehSymmetry}
h(-i,-j)=h(j,i),
\end{equation}
 which implies the symmetry condition
$$\overline{H(\xi,\eta)}=H(\eta,\xi).$$ Put
$$s(\xi,\eta)=|H(\xi,\eta)|^2=\sum_{i\in \ZZ}\sum_{j\in \ZZ}s_{ij}\xi^i\eta^j\;\;\;(s_{ij}\in \RR).$$
Then  $s:S^1\times S^1\rightarrow \RR$ is a kernel
of pure Fourier type satisfying Assumption~\ref{Assumption:Kernel}.
In particular, $\norm{s}_{L^1(C\times C)}=\norm{h}_{L^2(\ZZ\times \ZZ)}^2>0$.

\subsubsection{Filtered Wigner matrices (definition)}
For $i,j\in \{1,\ldots,N\}$ set
$$X_{ij}^{(N)}=\sum_{k=1}^N\sum_{\ell=1}^N Y_{k\ell}h(i-k,\ell-j),
$$
thus defining an $N$-by-$N$ random matrix $X^{(N)}$ which in view of the symmetry
(\ref{equation:LittlehSymmetry}) is real symmetric. We call $X^{(N)}/\sqrt{N}$ a {\em filtered Wigner matrix},
with {\em filter} $h$. We think of $X^{(N)}/\sqrt{N}$
as the result of filtering the Wigner matrix $Y^{(N)}/\sqrt{N}$ by $h$.

\subsubsection{Local covariance structure of $X^{(N)}$}
From (\ref{equation:YCorrelation}) and (\ref{equation:LittlehSupport})
we deduce that
\begin{equation}\label{equation:FilteredCorrelation}
i,j,k,\ell\in Q_{K}^{(N)}\;\&\;\min(j-i,\ell-k)>K\Rightarrow EX_{ij}^{(N)}X_{k\ell}^{(N)}=
s_{i-k,\ell-j},
\end{equation}
after a straightforward (extremely tedious) calculation. This is the main point in verifying that the real symmetric random matrices $X^{(N)}$
satisfy Assumption~\ref{Assumption:Model} with respect to the kernel $s$. The remaining details needed to check
Assumption~\ref{Assumption:Model} are easy to supply.
Thus we can put filtered Wigner matrices into the framework of our model.
\subsubsection{The (NE+SE+SW+NW)-filtered Wigner matrix}
\label{subsubsection:CompassPoints}
Taking 
$$h=(1/2)\one_{\{(1,-1),(1,1),(-1,1),(-1,-1)\}}, \;\;\mbox{and hence}\;\;s(e^{i\theta_1},e^{i\theta_2})=4\cos^2\theta_1\cos^2\theta_2,$$
we get a precisely defined version of the model mentioned in the 
introduction, which presently we will analyze in detail.
(We call this the (NE+SE+SW+NW)-filtered Wigner matrix
because each above the diagonal 
entry in the matrix
can be viewed as the average
of the four entries in a standard Wigner matrix which are to its immediate
northeast, southeast, southwest and northwest.)
Note that with $\theta$ uniformly distributed in $[0,2\pi)$,
the law of $2\cos^2\theta$ has density $\propto \one_{0<x<2}/\sqrt{x(2-x)}$.

\subsubsection{Remark}
 In \cite{AZ} we handled (real) Wishart matrices 
$Z^TZ$ (and more general matrices formed
from matrices $Z$ with independent real but perhaps not i.i.d.\ entries) 
in terms of  band-type matrices $\left[\begin{array}{cc}
0&Z\\
Z^T&0\end{array}\right]$.
A similar trick in the present setting
puts ``filtered Wishart matrices'' into the framework of our model.
For the kernels arising in that connection, 
the associated partition $\III$ consists of two intervals. We do not
discuss the Wishart case further here.

\subsection{The measure $\mu_s$}\label{subsection:MuDef}
We make the last preparation to state our main results.
Let $s$ be a kernel satisfying Assumption~\ref{Assumption:Kernel}.
For each positive integer $N$,  let 
$$C^{(N)}=\{c^{(N)}_1,\dots,c^{(N)}_{N^2}\}\subset C$$ 
be the set of pairs $(x,\xi)\in C$ where $x\in [0,1)\cap \frac{1}{N}\ZZ$ and $\xi^N=1$.
Then the empirical distribution $\frac{1}{N^2}\sum_{c\in C^{(N)}}\delta_c$
tends weakly as $N\rightarrow\infty$ to the uniform probability measure $P$. Let $C^{(\infty)}$
be the union of the sets $C^{(N)}$.
Let 
$$\{\tilde{Y}_e\}_{e\subset C^{(\infty)}\,\mathrm{s.t.}\, |e|=1,2}$$ 
be an i.i.d.\ family of standard normal (mean $0$ and variance $1$) random variables. (Recall from \S\ref{subsubsection:PureFourierType} that $|S|$ denotes the cardinality of $S$.)
Let $\tilde{X}^{(N)}$ be the $N^2$-by-$N^2$ real symmetric random matrix
with entries 
$$\tilde{X}_{ij}^{(N)}=2^{\delta_{ij}/2}\sqrt{s(c_i^{(N)},c_j^{(N)})}\tilde{Y}_{\{c_i^{(N)},c_j^{(N)}\}}.
$$
Let $\tilde{\lambda}_1^{(N)}\leq \cdots \leq \tilde{\lambda}_{N^2}^{(N)}$
be the eigenvalues of $\tilde{X}^{(N)}/N$ and let $\tilde{L}^{(N)}=\frac{1}{N^2}\sum_{i=1}^{N^2} \delta_{\tilde{\lambda}^{(N)}_i}$ be the empirical distribution
of the eigenvalues.  
By \cite[Thm.\ 3.2]{AZ} the empirical distribution $\tilde{L}^{(N)}$
tends weakly in probability as $N\rightarrow\infty$ to a limit $\mu_s$ 
with bounded support.
(The strange-looking factor $2^{\delta_{ij}/2}$ in the definition of $\tilde{X}^{(N)}$ could be dropped without changing the limit of $L^{(N)}$.
More generally, within the theory of \cite{AZ}, one has many ways to construct models with limiting measure $\mu_s$. We made our concrete choice to simplify the  proof of Proposition~\ref{Proposition:NonFluke} below.)
In Section~\ref{section:Method} we will provide a combinatorial description of the moments
of $\mu_s$ and  write down algebro-integral equations
(which we call {\em color equations})
satisfied by the Stieltjes transform of $\mu_s$ and certain auxiliary functions.

The following are our main results. In these results we fix a kernel
$s$ satisfying Assumption~\ref{Assumption:Kernel}
and a family $[X^{(N)}]_{N=1}^\infty$ of random hermitian matrices satisfying
Assumption~\ref{Assumption:Model} with respect to $s$.
As defined in \S\ref{subsubsection:EmpiricalDistribution}, let $L^{(N)}$ be the empirical distribution of the eigenvalues of $X^{(N)}/\sqrt{N}$.
Let $\mu=\mu_s$ be the measure
associated to $s$ by the procedure of \S\ref{subsection:MuDef}.

\begin{Theorem}\label{Theorem:MainResult}
$L^{(N)}$ converges weakly in probability to $\mu$.
 \end{Theorem}

\begin{Theorem}\label{Theorem:MainResultBis}
The Stieltjes transform $S(\lambda)=\int\frac{\mu(dx)}{\lambda-x}$ 
is an algebraic function of $\lambda$, i.~e.,
there exists some not-identically-vanishing polynomial $F(X,Y)$ in two variables with complex
coefficients such that $F(\lambda,S(\lambda))$ vanishes for all complex numbers $\lambda$ not in the support of $\mu$.
\end{Theorem}
\noindent To prove Theorem~\ref{Theorem:MainResult}
we first prove ``convergence in moments''
 in Section \ref{section:MainLimit},
and then finish the proof  in Section
\ref{section:BigCompletion} by considering variances. To prove Theorem~\ref{Theorem:MainResultBis} we
first set up a general
method  for proving algebraicity in Section \ref{section:Soft},
and then finish the proof in 
Section \ref{section:SoftApplied} by analyzing the color equations.

\subsection{Remarks on algebraicity}
$\;$

(i) In the setting of Theorem~\ref{Theorem:MainResultBis}, the polynomial $F(X,Y)$ is not unique, but if we impose the further condition of irreducibility,
then $F(X,Y)$ is unique up to a nonzero constant factor. (See \cite[Subsec.\ 2.2, Chap. 8]{ahl}.) When $F(X,Y)$ is so specified we call it
the {\em irreducible polynomial} for $S(\lambda)$.

(ii)
The method of proof of Theorem~\ref{Theorem:MainResultBis} does not yield an explicit irreducible polynomial $F(X,Y)$
for $S(\lambda)$.
But the elimination of variables necessary to make $F(X,Y)$ explicit could in principle be carried out on a computer using Gr\"{o}bner basis techniques.
For an introduction to the latter
see the very accessible book \cite{CLO}.
Our result serves as a theoretical guarantee that $F(X,Y)$ can indeed be found. 
The extent to which such Gr\"{o}bner basis calculations could be streamlined enough to be useful in applications remains to be investigated.


(iii) In recent work Rao and Edelman \cite{ER} explain how knowledge of 
an algebraic equation for the Stieltjes transform
can be exploited for numerical computation of limiting spectra in various random matrix models.
Furthermore, in the same setting, they explain how to streamline the computation of free additive and multiplicative convolution.

(iv) The general method of Section \ref{section:Soft}, 
and specifically Theorem  \ref{Theorem:BasicFunctionalEquation},
apply
to the Stieltjes transforms of the limiting measures
arising from the model of \cite{AZ}
in the case of a finite color space,
yielding algebraicity in all those cases. They also apply
to the  equations for the limit of the Stieltjes tranform
in \cite{HLN} and the systems of equations
in \cite[
Corollaries 9.3.2 and 10.1.2]{Gir},
under appropriate hypotheses of ``finiteness of filter'', as
well as many of the pre-limit ``Kx cannonical equations''
in \cite{Gir1}, x$\in \{1,\ldots,30\}$,
and in particular to his K1 cannonical equations.

%
%

(v) 
Questions of algebraicity have attracted attention for some time in areas bordering on random matrix theory.
For example, one may ask if the Green function for  random walk on a  free group is algebraic. 
A general result asserting this algebraicity under quite mild hypotheses was proved in \cite{ao}. 
These algebraicity questions were recently 
revisited in \cite{GB}, and a wide variety of striking connections were discussed, e.g., a connection with the theory of context-free languages.
An analogous problem in free probability
is to determine the spectra of convolution operators on a free group (and thus
through the theory of asymptotic freeness to determine the limit measures
for random matrix models involving several independent Haar-distributed unitary matrices in the large $N$ limit).
In the paper \cite{voi92}, which was primarily devoted to laying out
the  operator-valued version of the notion of freeness,
an analysis of such spectra was presented
as a sample application. 

(vi) Algebraicity implies a rather strong regularity property for $\mu$ which is of theoretical interest independent of computational issues.  We state and prove a theorem immediately below to explain this point in detail.

\subsection{Setup for the regularity theorem}
We declare a real-valued function $h$ defined in a bounded open interval $(a,b)$ to be  of {\em rational beta type} under the following conditions:
\begin{itemize}
\item $h$ is real-analytic and nonnegative on $(a,b)$.
\item There exist positive rational numbers $c$ and $d$ such that both limits
$$\lim_{x\downarrow a}(x-a)^{1-c} h(x),\;\;\;
\lim_{x\uparrow b}(b-x)^{1-d}h(x)$$ exist and are positive.
\end{itemize}
Note that $h\in L^1(a,b)$ and that
$h$ can have only finitely many zeroes in $(a,b)$.
For example $x^{c-1}(1-x)^{d-1}$ is of rational beta type 
in $(0,1)$ for all positive rational numbers $c$ and $d$ 
(which explains the terminology).

Now let $\mu$ be a probability measure on
 the real line with compact support $K$
and algebraic
Stieltjes transform  $S(z)=\int \frac{\mu(dx)}{z-x}$.
(We take $z$ as our complex variable instead of $\lambda$ 
to emphasize that we are now concerned only with complex analysis, not RMT.)
Let $F(X,Y)$ be a not-identically-vanishing polynomial 
such that $F(z,S(z))=0$ for all $z\in \CC\setminus K$ and furthermore
the discriminant $D(X)$ of $F(X,Y)$ with respect to $Y$ is 
not-identically-vanishing.
(We may, for example, take $F(X,Y)$ to be the irreducible polynomial of $S(z)$, 
because in that case the discriminant $D(X)$ cannot vanish identically.
But there is no need to insist on irreducibility of $F(X,Y)$; 
only the condition $D(X)\neq 0$ turns out to be important.)
Let $n$ be the degree of $F(X,Y)$ in $Y$ and 
write $F(X,Y)=\sum_{i=0}^n F_i(X)Y^i$. 
Let $A$ be the (finite) set of complex zeroes of $F_n(X)D(X)$.

\begin{Theorem}\label{Theorem:MainResultTer}
Notation and assumptions are as above.
Let $I$ be a connected component of $\RR\setminus A$ and let $\mu\vert_I$ be the restriction of $\mu$ to  (the Borel subsets of) $I$.
If $I\setminus K$ is nonempty, then $\mu\vert_I$ vanishes identically.
If $I$ is bounded and $\mu\vert_I$ does not vanish identically, 
then $\mu\vert_I$ has density of rational beta type with respect to Lebesgue measure.
\end{Theorem}
\noindent We emphasize that this result is essentially folkloric.
We supply a proof only for lack of a suitable reference.
We remark that the result does not rule out the possibility of $\mu$-atoms at  points of  $A\cap \RR$.
First we need a lemma.
For each positive integer $N$ let $r_N$ be a choice of 
analytic $N^{th}$ root of $z$ defined in the open
 set $\CC\setminus\{-i t\mid t\in [0,\infty)\}$.
Throughout, we use $\Im z$ for the imaginary part of $z\in\CC$.

\begin{Lemma}
Choose $\delta>0$ so small that
$A\cap \RR=A\cap\{|\Im z|<\delta\}$.
Let 
$$W=\{|\Im z|<\delta\}\setminus
\bigcup_{a\in A\cap \RR}
\{a-i t\mid t\in [0,\infty)\}.
$$
Then: (i) There exists a unique analytic function $g$ defined in $W$ which agrees with $S$ on the set $\{0<\Im z<\delta\}$.
(ii) For each $a\in A\cap \RR$
there exists a positive integer $N$ and $\epsilon>0$ such that
$(z-a)^Ng(z)$ admits expansion
in the set $W\cap \{|z-a|<\epsilon\}$
as a convergent power series in $r_N(z-a)$.\end{Lemma}
\proof Let $\hat{F}(X,Y)$ be the irreducible polynomial of $S(z)$, which divides $F(X,Y)$.
Let $\hat{n}$ be the degree of $\hat{F}(X,Y)$
in $Y$ and write $\hat{F}(X,Y)=\sum_{i=0}^{\hat{n}}\hat{F}_i(X)Y^i$.
Let $\hat{D}(X)$ be the discriminant of $\hat{F}(X,Y)$ with respect to $Y$.
Let $\hat{A}$ be the set of complex zeroes of $\hat{F}_{\hat{n}}(X)\hat{D}(X)$.
Then we have $\hat{A}\subset A$.
Let $\mathbf{f}$  be the collection of all pairs $(f,U)$ consisting of a open subset $U\subset \CC$ and an analytic function $f$ defined in $U$ such that $\hat{F}(z,f(z))=0$ for all $z\in U$. 
Then $\mathbf{f}$ is the global analytic and algebraic function naturally associated to $\hat{F}(X,Y)$.
(Here up to minor changes of notation we follow the definitions of
 \cite[Secs.\  1,2, Chap.\ 8]{ahl}.) The main point of this construction for us is that
the pair $(S,\CC\setminus K)$ belongs to $\mathbf{f}$.
Every branch of $\mathbf{f}$ at a point in $\CC\setminus \hat{A}$ 
admits analytic continuation along any arc 
not passing through the set $\hat{A}$.
(See \cite[Subsec. 2.2, Pg. 303, Chap. 8]{ahl}.)
Since $U$ and $\{0<\Im z<\delta\}$ are 
simply connected and disjoint from $\hat{A}$,
statement (i) follows by the Monodromy Theorem \cite[Theorem 2, Subsec.\ 1.6,
Chap. \ 8]{ahl}.
Statement (ii) follows from the study of 
singularities of algebraic functions in
\cite[Subsec.\ 2.3, Chap.\ 8]{ahl}. 
\qed

\proof[Proof of the theorem]
Recall the inversion formula
$$
\int \phi(x)\mu(dx)=\frac{1}{\pi}\lim_{y\downarrow 0}\int \Im S(x+i y)\phi(x)dx
$$
holding for all bounded continuous functions $\phi$ defined on the real line.  With $g$ as in part (i) of the previous lemma, let $h=\frac{1}{\pi}\Im g\vert_{\RR\setminus A}$.
From the inversion formula it follows that $\mu\vert_I$ has density $h\vert_I$ with respect to Lebesgue measure. If $I\setminus K$ is nonempty (this is in particular the case if $I$ is unbounded), then
$h\vert_I$ vanishes identically in some (nonempty) open subinterval of $I$,
and hence (since real analytic) vanishes identically. Otherwise $I$ is bounded,
and if $h\vert_I$ does not vanish identically,
then $h\vert_I$ is of rational beta type by part (ii) of the preceding lemma
along with the fact that $h\vert_I\in L^1(I)$.
\qed

\section{The moments and Stieltjes transform of $\mu$}
\label{section:Method}
We fix a kernel $s$ satisfying Assumption~\ref{Assumption:Kernel}
and put $\mu=\mu_s$. We provide a detailed description of the moments and Stieltjes transform of $\mu$. We also introduce combinatorial tools needed throughout the paper.

\subsection{Quick review of key combinatorial notions}
\label{sec-2.2}

\subsubsection{Graphs}
For us a {\em graph} $G=(V,E)$ consists by definition of a finite set
$V$ of {\em vertices} and a set $E$
of {\em edges}, where every element of $E$
is a subset of $V$ of cardinality $1$ or $2$.
In other words, we are dealing here with graphs 
(i) which have finitely many vertices, (ii) which have unoriented edges,
(iii) in which a vertex may be joined to itself by an edge,
but (iv) in which no two vertices may be joined by more than one edge.
A graph $G$ is a {\em tree} if $G$ is connected
and $|V|=|E|+1$. 
 We emphasize that every edge of a tree  joins two distinct vertices---it is not allowed for a
vertex of a tree to be joined to itself by an edge.

\subsubsection{Set partitions}
We say that a set $\pi\subset 2^{\{1,\dots,k\}}$ is a {\em set-partition} of $k$
if $\pi$ is a disjoint
family of nonempty sets with union equal to $\{1,\dots,k\}$. The elements of $\pi$ 
are called the {\em parts} of $\pi$.
For each $i\in \{1,\dots,k\}$, let $\pi(i)$
be the part of $\pi$ to which $i$ belongs.
For convenience we extend $i\mapsto \pi(i)$
to a periodic function on $\ZZ$ by the rule $\pi(i)=\pi(i+k)$.

\subsubsection{The graph associated to a set partition}
To each set partition $\pi$ of $k$ there is canonically associated a graph $G_\pi=(V_\pi,E_\pi)$,
where $V_\pi=\pi$ and 
$$E_\pi=\{\{\pi(i),\pi(i+1)\}\mid i=1,\dots,k\}.$$
By construction $G_\pi$ comes canonically equipped with a walk, namely $$\pi(1),\dots,\pi(k),\pi(k+1)=\pi(1),$$
 whence in particular it follows that $G_\pi$ is connected. 

\subsubsection{Wigner set partitions}
We say that a set partition $\pi$ of $k$ is a {\em Wigner set partition}
if the corresponding graph $G_\pi$
has $k/2+1$ vertices and $k/2$ edges,
in which case $G_\pi$,
since connected, is a tree.
We denote the set of such $\pi$ by $\WWW_k$.
For $k$ odd the set $\WWW_k$ is empty.
For $k=2\ell$ the set $\WWW_k$ is canonically in bijection with the set of rooted planar trees
with $\ell+1$ nodes
and hence, as is well-known \cite{stanley},
the cardinality
of $\WWW_k$ is the {\em Catalan number} 
$\frac{1}{\ell+1}\left(\begin{subarray}{c}
2\ell\\
\ell
\end{subarray}\right)$. 
\begin{Lemma}\label{Lemma:CrucialInsight}
Fix $\pi\in \WWW_k$. (i) For each $i\in \{1,\dots,k\}$ we have
$\pi(i)\neq \pi(i+1)$.  \linebreak (ii) For each $e\in E_\pi$ 
the equation $e=\{\pi(i),\pi(i+1)\}$ has exactly two solutions \linebreak $i\in \{1,\dots,k\}$,
say $i_1$ and $i_2$, and moreover $\pi(i_1)\neq \pi(i_2)$. (iii) For each \linebreak $i\in \{1,\dots,k\}$ we have $\{\pi(i),\pi(i+1)\}=\{\pi(j-1),\pi(j)\}$,
where $j$ is the least of the integers $\ell>i$ such that $\pi(i)=\pi(\ell)$.
\end{Lemma}
\proof The lemma formalizes facts about the tree $G_\pi$ and the canonical walk on it which are clear from a graph-theoretic point of view. 
(i)
No edge of $G_\pi$ connects a vertex to itself.
(ii) The canonical walk on $G_\pi$ visits each edge of $G_\pi$
exactly twice. More precisely, the canonical walk traverses each edge of $G_\pi$ exactly once in each direction.  (iii) The canonical walk on $G_\pi$ extended by periodicity
returns to a given vertex on the same edge by which it departed.
\qed

\subsubsection{Tree integrals}
Let $\{\kappa_A\}$ be an i.i.d.\ family of $C$-valued uniformly distributed 
 random variables indexed by finite nonempty sets of positive integers. 
Expectations with respect to these variables are denoted by $\E$.
For each $\pi\in \WWW_k$ we define a bounded random variable by the formula
$$M_\pi=\prod_{\{A,B\}\in E_\pi}s(\kappa_A,\kappa_B),$$
which is  well-defined on account of the symmetry $s(c,c')=s(c',c)$.
We call the expectation $\E M_\pi$ a {\em tree integral}.

\begin{Proposition}[Combinatorial description of the moments of $\mu$]
\label{Proposition:MuMoments}
We have
\begin{equation}\label{equation:MuMomentsBis}
\langle \mu,x^k\rangle=\sum_{\pi\in \WWW_k} \E M_\pi
\end{equation}
for every integer $k>0$. 
\end{Proposition}
\noindent The proposition is essentially just a special case of \cite[Thm.\ 3.2]{AZ},
but a fair amount of explanation is needed because the set up in this paper 
is (superficially) incompatible with that of \cite{AZ}---here we emphasize set partitions,
whereas in \cite{AZ} we emphasized ``words'' and ``spelling''.
We can immediately reduce
the proposition to the following technical lemma.
The lemma is slightly more detailed than needed for the proof of the proposition---part (ii) will be needed for the derivation of algebro-integral equations for the Stieltjes transform of $\mu$.

\begin{Lemma}\label{Lemma:MuMoments}
Put $A=2\norm{s}_{L^\infty(C\times C)}^{1/2}$.
(i) There exists a unique system \linebreak $\{\Phi_n,\Psi_n\}_{n\in\NN}$
of functions in $L^\infty(C)$ such that
\begin{equation}\label{equation:FirstPhiPsi}
\Psi_n(c)=\int s(c,c')\Phi_n(c')P(dc')
\end{equation}
holds $P$-a.e.\ for every $n$ and 
there holds an identity
\begin{equation}\label{equation:SecondPhiPsi}
\sum_{n=1}^\infty \Phi_nt^n=t\left(1-t\sum_{n=1}^\infty \Psi_nt^n\right)^{-1}
\end{equation}
of formal power series in $t$ with coefficients in $L^\infty(C)$.  
(ii) The bounds
\begin{equation}\label{equation:PreMcGuffin}
0\leq \Phi_{n}\leq A^{n-1},\;\;\;\;
0\leq \Psi_{n}\leq A^{n+1}/4
\end{equation}
hold $P$-a.e.\ for every  $n$.
(iii) The formula 
\begin{equation}\label{equation:McGuffin}
\langle P,\Phi_{k+1}\rangle=\sum_{\pi\in \WWW_{k}}\E M_\pi
\end{equation}
holds for every integer $k>0$.
\end{Lemma}
\proof[Proof of the proposition granting the lemma]
According \cite[Lemma 3.2]{AZ} (taking there $\sigma=P$, $D=0$, and $s^{(2)}=s$),
there exists a unique probability 
measure on the real line with $k^{th}$ moment $\langle P,\Phi_{k+1}\rangle$
for every $k\geq 0$. That measure according to \cite[Thm.\ 3.2]{AZ}
is none other than $\mu$.
\qed

\proof[Plan for the proof of the lemma]
Part (i) of the lemma is nothing but an inductive definition
of $\Phi_n$ and $\Psi_n$.
So only parts (ii,iii) require proof.
In principle, part (iii) follows from \cite[Lemmas 6.3 and 6.4]{AZ},
but because of the incompatibility of set-ups noted above, those lemmas 
cannot be directly applied here---some amplification is needed.
Also part (ii) is most easily explained from the point of view of
\cite[loc.\ cit.]{AZ}
So, after recalling in \S\ref{subsection:Verbal} 
the needed background from \cite{AZ},
we lightly sketch  in \S\ref{subsection:LightSketch} a proof of parts (ii,iii) of the lemma.

\subsection{The ``verbal'' approach}\label{subsection:Verbal}
We briefly recall the  point of view of \cite{AZ} and 
compare it to the present one. The material reviewed 
here  will be used in a substantial way only in Section \ref{section:Method},
not in later sections of the paper.

\subsubsection{Words}
We fix a set of {\em letters} and define a {\em word}
to be a finite nonempty sequence $w=\alpha_1\cdots\alpha_k$ 
of letters.  We say that words $w=\alpha_1\cdots \alpha_k$
and $x=\beta_1\cdots \beta_\ell$ are {\em equivalent}
if $k=\ell$ (the words are the same length) and there exists
a one-to-one-correspondence $\varphi:\{\alpha_i\}\rightarrow\{\beta_j\}$
such that $\varphi(\alpha_i)=\beta_i$ for $i=1,\dots k$ 
(each word codes to the other under a simple substitution cipher).
Each word $w=\alpha_1\cdots\alpha_k$ of length $k$ gives rise naturally to a
set partition of $k$, namely the set of equivalence classes for the 
equivalence relation $i\sim j\Leftrightarrow \alpha_i=\alpha_j$.
Two words are equivalent if and only if they have the same
length and give rise to the same set partition. The upshot is that
speaking of equivalence classes of words is 
equivalent to speaking of set partitions.

\subsubsection{Wigner words}
Let $w$ be a word of at least two letters with 
same first and last letter, and let $w'$ be the word
obtained by dropping the last letter of $w$.
The word $w$ is a Wigner word in the sense of \cite{AZ}
if and only if the set partition associated to $w'$ is a Wigner set partition
in the sense of this paper. In \cite{AZ} we also declared every
one-letter word to be a Wigner word.  The Wigner words have a 
simple inductive characterization \cite[Prop.\ 4.5 and \S4.7]{AZ}  which is not so convenient to state
in the set partition language.
 To wit, a word $w$ is a Wigner word if and only if the
 following conditions hold:
\begin{itemize}
\item The first and last letters of $w$ are the same.
\item No letter appears twice in a row in $w$.
\item Let $\alpha$ be the first letter of $w$. 
Write $w=\alpha w_1\alpha \cdots \alpha w_r\alpha$, where $\alpha$ does not
appear in any of the words $w_i$. Then each word $w_i$ is a Wigner word,
 and moreover for $i\neq j$ the words $w_i$ and $w_j$ have no letters in common.
(If $r=0$ then $w$ consists of a single letter and is by definition
 a Wigner word.)
\end{itemize}

\subsubsection{``Verbal'' description of tree integrals}
Let $\{\kappa_\alpha\}$ be a letter-indexed i.i.d.\ family of
 uniformly distributed
$C$-valued random variables. Given a Wigner word $w$, we define
a random variable $M_w$ inductively by the following procedure:
\begin{itemize}
\item Writing $w=\alpha w_1\alpha \cdots \alpha w_r \alpha$
as in the inductive characterization of Wigner words, 
and letting $\beta_i$ be the initial
letter of $w_i$ for $i=1,\dots,r$, we set
$M_w=\prod_{i=1}^r s(\kappa_\alpha,\kappa_{\beta_i})M_{w_i}$. 
(When $w$ is one letter long, we put $M_w=1$.)
\end{itemize}
The formula (\ref{equation:McGuffin}) claimed in 
Lemma~\ref{Lemma:MuMoments} can 
be rewritten
\begin{equation}\label{equation:MomentRewrite}
\langle P,\Phi_{k+1}\rangle=\sum_{w\in W_{k+1}} \E M_w
\end{equation}
where the sum is extended over a set  of representatives $W_{k+1}$
for equivalence classes of Wigner words
of length $k+1$. Further and crucially,
notation as above in the inductive characterization of the random variables $M_w$, we have a relation
\begin{equation}\label{equation:TreeIntegralStructure}
\E (M_w\vert \kappa_\alpha)=
\prod_{i=1}^r \E(s(\kappa_\alpha,\kappa_{\beta_i})
\E( M_{w_i}\vert \kappa_{\beta_i})\vert \kappa_\alpha)
\end{equation}
among conditional expectations.

\subsection{Proof of 
Lemma~\ref{Lemma:MuMoments}(ii,iii)}\label{subsection:LightSketch}
It is enough to construct an example
of a system $\{\Phi_n,\Psi_n\}$ in $L^\infty(C)$ satisfying
(\ref{equation:FirstPhiPsi},\ref{equation:SecondPhiPsi},\ref{equation:PreMcGuffin},\ref{equation:McGuffin}),
and to do so we follow the path of the proofs of \cite[Lemmas 6.3 and 6.4]{AZ}. 
Fix a letter $\alpha$. 
We may suppose that every word belonging
to the set of representatives $W_{k+1}$ figuring in 
(\ref{equation:MomentRewrite})
begins with $\alpha$. There exist for each integer $k\geq 0$ 
well-defined $\Phi_{k+1},\Psi_{k+1}\in L^\infty(C)$ such that
$$\Phi_{k+1}(\kappa_\alpha)=\sum_{w\in W_{k+1}}\E (M_w\vert \kappa_\alpha),\;\;\;
\Psi_{k+1}(\kappa_{\alpha'})=\sum_{w\in W_{k+1}}
\E(s(\kappa_{\alpha'},\kappa_\alpha)M_w\vert \kappa_{\alpha'}),$$
where $\alpha'$ is a letter not appearing in any of the words belonging to the set $W_{k+1}$.
The system $\{\Phi_n,\Psi_n\}$
 has property 
(\ref{equation:FirstPhiPsi}) by construction,
and has property (\ref{equation:McGuffin}) since
 the latter is equivalent to (\ref{equation:MomentRewrite}).
Since 
$$|\WWW_k|=|W_{k+1}|\leq 2^k,\;\;\;
0\leq M_w\leq \norm{s}_{L^\infty(C\times C)}^{k/2}\;\mbox{for}\;w\in W_{k+1},$$ 
for all integers $k\geq 0$, the system $\{\Phi_n,\Psi_n\}$
has property (\ref{equation:PreMcGuffin}).
Finally,  from (\ref{equation:TreeIntegralStructure})
and the inductive characterization of Wigner words, we obtain identities
$$\Phi_{k+1}=
\sum_{r=0}^\infty\sum_{
\begin{subarray}{c}(\ell_1,\dots,\ell_r)\in \NN^r\\
\sum_{i=1}^r (\ell_i+1)=k
\end{subarray}}\prod_{i=1}^r\Psi_{\ell_i}$$
 in $L^\infty(C)$ for all integers $k\geq 0$ 
which together imply that the system
$\{\Phi_n,\Psi_n\}$ has property (\ref{equation:SecondPhiPsi}).
The proofs of Lemma~\ref{Lemma:MuMoments} and 
Proposition~\ref{Proposition:MuMoments} are now complete.
\qed

\subsection{The  color equations}
We continue in the setting  of Proposition~\ref{Proposition:MuMoments}.
\subsubsection{Nice functions}
We say that a complex-valued function $f$ on color space
is {\em nice} (with respect to the kernel $s$ and associated partition $\III$ of $[0,1]$) if $f$ has a Fourier expansion 
$$f(c)=\sum_{i\in \ZZ} f_i(x)\xi^i\;\;(c=(x,\xi)\in C)$$
where all but finitely many of the coefficient functions $f_i:[0,1]\rightarrow \CC$
vanish identically, and every coefficient function $f_i$
is constant on every interval of the partition $\III$.
It is not difficult to see that $\Phi_n$ and $\Psi_n$ 
can be ``corrected'' on a set of $P$-measure
zero in a unique way to become nice.
So we may and we will assume hereafter without
any loss of generality
that the functions $\Phi_n$ and $\Psi_n$ are nice and that all the relations asserted in Proposition~\ref{Proposition:MuMoments} to hold $P$-a.e.\ in fact hold without exception.

\subsubsection{An auxiliary function}
For all  complex numbers $|\lambda|>A$ and $c\in C$ put
$$\Psi(c,\lambda)=\sum_{n=1}^\infty \Psi_n(c)\lambda^{-n},$$
defining a function depending holomorphically
on $\lambda$ and satisfying estimates
\begin{equation}\label{equation:PsiLambdaEst}
|\Psi(c,\lambda)|\leq \frac{1}{4}\frac{A^2}{|\lambda|-A},\;\;\;
|\lambda|>2A\Rightarrow |\Psi(c,\lambda)|<\frac{A^2}{2|\lambda|}\leq \frac{A}{4}
\end{equation}
uniform in $c$.

\subsubsection{The equations}
From the system of equations described in Proposition~\ref{Proposition:MuMoments}
we now deduce by the substitution $t=1/\lambda$ and application of dominated convergence the relations
\begin{equation}\label{equation:Color}
\int \frac{s(c,c')P(dc')}{\lambda-\Psi(c',\lambda)}=\Psi(c,\lambda),\;\;\;
\int\frac{P(dc)}{\lambda-\Psi(c,\lambda)}=S(\lambda)=\int \frac{\mu(dx)}{\lambda-x}
\end{equation}
which hold for all $c\in C$ and $|\lambda|>2A$.
We call the relations (\ref{equation:Color}) the {\em  color equations}.
Equations of this sort have appeared already in other contexts, see \cite{BK99},
\cite{Gir}, \cite[eq. 2.2]{HLN}, \cite{khorunzhy1}.
 
\begin{Proposition}\label{Proposition:NonFluke}
In the setting of the  color equations, assume further that
for some nice nonnegative function $f$ on color space
$$s(c,c')=f(c)f(c'),\;\;\;\norm{f}_{L^\infty(C)}= A/2,\;\;\;\norm{f}_{L^1(C)}=1.$$
Let $\mu_f$ be the law of $f$ viewed as a random variable
on $C$ under $P$.
Let $S_f(\lambda)$ be the Stieltjes transform of $\mu_f$.
Then: (i) There exists a function $w(\lambda)$ defined and holomorphic
for $|\lambda|\gg 0$ such that 
\begin{equation}\label{equation:PreRankOneMaster}\lim_{|\lambda|\rightarrow\infty} w(\lambda)\lambda=1,
\end{equation}
\begin{equation}\label{equation:RankOneMaster}
\lambda S(\lambda)=1+w(\lambda)^2=\frac{\lambda}{w(\lambda)}S_f\left(
\frac{\lambda}{w(\lambda)}\right)
\end{equation}
for $|\lambda|\gg 0$. (We emphasize that we do mean $S$ on the LHS and $S_f$ on the RHS.) \linebreak
(ii) The measure $\mu$ is the free multiplicative convolution
of $\mu_f$  with the semicircle law of mean $0$ and variance $1$.
\end{Proposition}
\proof (i) Put
\begin{equation}\label{equation:TopDisplay}
w(\lambda)=\int \frac{f(c)P(dc)}{\lambda-\Psi(c,\lambda)}=\int \Psi(c,\lambda)P(dc),
\end{equation}
thus defining a holomorphic function in the domain $|\lambda|>2A$ such that
$$w(\lambda)=\frac{1}{\lambda}+O\left(\frac{1}{|\lambda|^2}\right),\;\;\;\Psi(c,\lambda)=w(\lambda)f(c).$$
By definition
$$S_f(\lambda)=\int\frac{P(dc)}{\lambda-f(c)},$$
which is a function holomorphic in the domain $|\lambda|>A/2$.
For $|\lambda|\gg 1$ we have by (\ref{equation:Color})
and the first equality in (\ref{equation:TopDisplay}) that
$$
\begin{array}{rcl}
\lambda S(\lambda)&=&
\displaystyle \int\frac{\lambda P(dc)}{\lambda-w(\lambda)f(c)}=\frac{\lambda}{w(\lambda)}S_f\left(\frac{\lambda}{w(\lambda)}\right),\\\\
\lambda S(\lambda)-1
&=&
\displaystyle\int\frac{w(\lambda)f(c) P(dc)}{\lambda-w(\lambda)f(c)}=w(\lambda)^2,
\end{array}
$$
which proves the result. 

(ii) We return to the setting of \S\ref{subsection:MuDef}.
Let $W^{(N)}$ be the $N^2$-by-$N^2$ matrix with entries $$W_{ij}^{(N)}=2^{\delta_{ij}/2}
\tilde{Y}_{\{c_i^{(N)},c_j^{(N)}\}}/N.$$
Note that $W^{(N)}$ belongs to the GOE (the factor
$2^{\delta_{ij}/2}$ is needed for orthogonal invariance).
Let $\Lambda^{(N)}$ be the $N^2$-by-$N^2$ deterministic diagonal matrix
with diagonal entries 
$$\Lambda_{ii}^{(N)}=\sqrt{f(c_i^{(N)})}.$$
Then we have
$$\tilde{X}^{(N)}/N=\Lambda^{(N)}W^{(N)}\Lambda^{(N)}.$$
Furthermore, as $N\rightarrow\infty$,
the empirical distribution of eigenvalues of $W^{(N)}$ (resp., $(\Lambda^{(N)})^2$)
tends to the semicircle law (resp., $\mu_f$).  Finally, since $W^{(N)}$
and $\Lambda^{(N)}$ are asymptotically freely independent,
see \cite[Corollary 4.3.6]{HP} and the discussion on page
157 there concerning the extension from the unitary to orthogonal
case, 
$\mu$ has the claimed form of free multiplicative convolution.
\qed

\subsection{Analysis of (NE+SE+SW+NW)-filtered Wigner matrix}\label{section:Example}
We consider the 
setup of \S\ref{sec-motivation} in the special case
mentioned in \S\ref{subsubsection:CompassPoints}.
We are thus considering the model mentioned in the introduction.
Proposition~\ref{Proposition:NonFluke} applies,
with $\mu_f$ equal to the law of $2\cos^2\theta=1+\cos 2\theta$
with $\theta$ distributed uniformly in $[0,2\pi)$.
Integrating, we get
$$S_f(\lambda)=\frac{1}{2\pi} \int_0^{2\pi} \frac{d\theta}{\lambda-1-\cos 2\theta}=
\frac{1}{\sqrt{\lambda(\lambda-2)}}\,.$$
It follows from (\ref{equation:RankOneMaster})
that 
$$(1+w^2)^2=\frac{\lambda}{\lambda-2 w},$$
and hence (after taking out an irrelevant factor of $w$)
$$
2w^4-\lambda w^3+4 w^2-2\lambda w+2=0.
$$
In turn, after forming the resultant
of lefthand side above and $1+w^2-\lambda S$ with respect to $w$ 
and  taking out irrelevant factors,
we get the  equation
\begin{equation}
\label{eq-140706}
4\lambda^2S^4-\lambda^3S^3-\lambda^2S^2+\lambda S+1=0.
\end{equation}
Since the equation is quartic, it can 
in principle be solved explicitly by root extractions.
We omit the details, which are available from the authors. 
Fortunately explicit formulas for $S(\lambda)$ are not needed to find out the key features of $\mu$.
We can use Theorem~\ref{Theorem:MainResultTer} and its proof to get information without having to calculate too much.
The discriminant of the
left side of (\ref{eq-140706}) is
$$-16\lambda^6(8 \lambda^4+107\lambda^2-1024).$$
The only nonzero real roots of the discriminant are
\begin{equation}\label{equation:EndpointSurds}
\pm \frac{1}{4}\sqrt{-107+51\sqrt{17}}\;\;\;
 (\mbox{approximately}\;\pm 2.5406).
 \end{equation}
Since $\mu$ is symmetric and cannot be supported on finitely many points (equation (\ref{eq-140706}) cannot be solved by any rational function of $\lambda$), 
Theorem~\ref{Theorem:MainResultTer} leaves no choice but that the support of $\mu$  be the interval with  endpoints specified in (\ref{equation:EndpointSurds}).
To finish our discussion let us verify that that $d\mu/dx$ has a spike at the origin proportional to $1/\sqrt{|x|}$.
Let $\sqrt{\lambda}$ be the unique square root of $\lambda$
in $\CC\setminus -i[0,\infty)$ which is positive along the positive real axis. For small enough $\epsilon>0$ 
there are exactly four analytic solutions of
(\ref{eq-140706}) 
in the slit disk $\{|\lambda|<\epsilon\}\setminus -i[0,\infty)$, 
and these are of the form $i^\nu(\sqrt{i}/\sqrt{2\lambda})+O_{\lambda\rightarrow 0}(1)$ for $\nu=0,1,2,3$. This is verified by using the fact that such solutions have in any case for $\epsilon>0$ sufficiently small
 an expansion on $\{|\lambda|<\epsilon\}\setminus -i[0,\infty)$ in fractional (possibly negative) powers of $\lambda$, and 
then one calculates the expansions by the method of undetermined coefficients.
The only such solution with positive imaginary part
on $(-\epsilon_1,0)\cup(0,\epsilon_1)$ for $\epsilon_1>0$ sufficiently small
corresponds to $\nu=0$; this gives the spike.

\section{The main limit calculation}\label{section:MainLimit}
To the end of proving Theorem~\ref{Theorem:MainResult},
we first prove the following result.

\begin{Proposition}\label{Proposition:MainLimit}
Let Assumption~\ref{Assumption:Model} hold.
For each positive integer $k$,
$$\lim_{N\rightarrow\infty}
N^{-k/2-1}E\trace ((X^{(N)})^k)=
\sum_{\pi\in \WWW_k}\E M_\pi.$$
\end{Proposition}
\noindent The proof requires some preparation and will not be completed until \S \ref{subsection:MainLimit}.

\subsection{Notation, terminology and strategy}
\subsubsection{$(N,k)$-words}
Let $N$ and $k$ be positive integers.
An {\em $(N,k)$-word} $\ibold$ is by definition a function 
$$\ibold:\{1,\dots,k\}\rightarrow \{1,\dots,N\}.$$
To each $(N,k)$-word we attach a random variable
$$X_{\ibold}^{(N)}=\prod_{\alpha=1}^k X_{\ibold(\alpha),\ibold(\eta_k(\alpha))}^{(N)}$$
where
$$\eta_k=(12\cdots k)\in S_k.$$
(Here and below we employ cycle notation for permutations.)
We have
$$
\sum_{\ibold}E X_\ibold^{(N)}=E\trace ((X^{(N)})^k),
$$
where the sum on the left is extended over $(N,k)$-words $\ibold$.

\subsubsection{The set partition associated to an $(N,k)$-word}\label{subsubsection:BlobDef}
Given an $(N,k)$-word $\ibold$,
consider the graph $G^K_\ibold=(V^K_\ibold,E^K_\ibold)$,
where $V^K_\ibold=\{1,\dots,k\}$
and 
$$E^K_\ibold=\{\{\alpha,\beta\}\subset
\{1,\dots,k\}\mid |\ibold(\alpha)-\ibold(\beta)|\leq K\}.$$
Here and below $K$ is the  constant figuring in Assumption~\ref{Assumption:Model}(II).
We define $\pi_\ibold$ to be the set partition of $k$
the parts of which are the equivalence classes
for  the  relation 
``$\alpha$ and $\beta$ belong to the same
connected component of $G^K_\ibold$''.
(Although the set partition $\pi_\ibold$ depends on $K$, we suppress reference to $K$ in the notation.)
We have 
$$|\ibold(\alpha)-\ibold(\beta)|\leq K\Rightarrow \pi_\ibold(\alpha)=\pi_\ibold(\beta),\;\;\; |\ibold(\alpha)-\ibold(\beta)|>Kk\Rightarrow \pi_\ibold(\alpha)\neq \pi_\ibold(\beta),
$$
for all $\alpha,\beta\in \{1,\dots,k\}$. 

\subsubsection{Distinguished $(N,k)$-words}
Let $\ibold$ be an $(N,k)$-word.  
Consider the following conditions:
\begin{enumerate}
\item $\pi_\ibold\in \WWW_k$.
\item For all $A\in \pi_\ibold$ we have $\ibold(\min A)\in Q^{(N)}_{(k+1)K}$.
\item For all distinct $A,B\in \pi_\ibold$
we have $|\ibold(\min A)-\ibold(\min B)|>5kK$.
\item $\displaystyle\max_{A\in \pi_\ibold}
\max_{\alpha,\beta\in A}|\ibold(\alpha)-\ibold(\beta)|\leq Kk$.
\item $\displaystyle\min_{\alpha=1}^k
|\ibold(\alpha)-\ibold(\eta_k(\alpha))|>3Kk$.
\item For $\alpha=1,\dots,k$ we have $\ibold(\alpha)\in Q^{(N)}_{K}$.
\item For all $A\in \pi$ and $\alpha,\beta\in A$,
the numbers $\frac{\ibold(\alpha)}{N}$ and
$\frac{\ibold(\beta)}{N}$ belong to the same
interval of the partition $\III$.
\end{enumerate}
If conditions (I,II,III) hold we say that $\ibold$ is {\em distinguished}, in which case $\ibold$ automatically also satisfies conditions (IV,V,VI,VII).

\subsubsection{Strategy}
We will show that the only nonnegligible contributions to
\begin{equation}\label{equation:Promise}
\lim_{N\rightarrow\infty}
N^{-k/2-1}\sum_{\mbox{\scriptsize
$\ibold$: $(N,k)$-word}} EX_\ibold^{(N)}=
\lim_{N\rightarrow\infty}
N^{-k/2-1}E\trace((X^{(N)})^k)
\end{equation}
come from distinguished $(N,k)$-words.
Then we will evaluate 
$EX^{(N)}_\ibold$ for distinguished $\ibold$.
Finally, we will calculate the limit on the left with the summation restricted to distinguished $\ibold$.

\subsection{Negligibility of nondistinguished $(N,k)$-words}
\label{subsection:Suite}
\begin{Lemma}\label{Lemma:CrudeCounting}
Let $\pi$ be a set partition of $k$. 
There exists $C_\pi>0$
such that for every positive integer $N$
the sum $\sum_\ibold E|X^{(N)}_\ibold|$
extended over  $(N,k)$-words $\ibold$
such that $\pi=\pi_\ibold$ does not
exceed $C_{\pi} N^{|\pi|}$. 
\end{Lemma}
\proof  By Assumption~\ref{Assumption:Model}(Ib)
and the H\"{o}lder inequality, it suffices simply to estimate the number of $(N,k)$-words
such that $\pi=\pi_\ibold$. A crude estimate
of the latter  is $(1+2Kk)^{k-|\pi|}N^{|\pi|}$. \qed

\begin{Lemma}\label{Lemma:WignerBound}
Let $\pi$ be a set partition of $k$.
Let $\ibold$ be an $(N,k)$-word such that
$\pi=\pi_\ibold$.
If $|\pi|\geq k/2+1$ and $EX_\ibold^{(N)}\neq 0$,
then $\pi\in \WWW_k$ (and hence $|\pi|=k/2+1$).
\end{Lemma}
\proof Let $G_\pi=(V_\pi,E_\pi)$ be the graph associated to $\pi$.   For each $\alpha\in \{1,\dots,k\}$
put $e(\alpha)=\{\pi(\alpha),\pi(\alpha+1)\}\in E_\pi$.
Now fix $e\in E_\pi$. It is enough to show
that $e=e(\alpha)$ for at least two $\alpha\in \{1,\dots,k\}$, for then, since $G_\pi$ is connected, we have 
$$k/2+1\leq |\pi|=|V_\pi|\leq |E_\pi|+1\leq k/2+1,$$
and hence $\pi$ is a Wigner set partition. To derive a contradiction, suppose rather that $e=e(\alpha)$
for unique $\alpha\in \{1,\dots,k\}$.
For every $\gamma\in\{1,\dots,k\}$
let $i(\gamma)\leq j(\gamma)$ be the integers
$\ibold(\gamma)$ and $\ibold(\eta_k(\gamma))$ rearranged.
Then for every $\beta\in \{1,\dots,k\}\setminus\{\alpha\}$ we have 
$\max(|i(\alpha)-i(\beta)|,|j(\alpha)-j(\beta)|)>K$, for otherwise $e(\alpha)=e(\beta)$. 
It follows by Assumption~\ref{Assumption:Model}(IIb)
that the random variable $X^{(N)}_{\ibold(\alpha),\ibold(\eta_k(\alpha))}$
is independent of the rest of the random variables
appearing in the product $X^{(N)}_\ibold$,
and hence $EX_\ibold^{(N)}=0$ by Assumption~\ref{Assumption:Model}(Ia).
This contradiction proves the lemma. \qed
 
\begin{Lemma}\label{Lemma:CrudeCountingBis}
Let $\pi$ be  a Wigner set partition.
There exists $C'_\pi>0$ such that for every  positive integer $N$ the
sum $\sum_\ibold E|X^{(N)}_\ibold|$
extended over $(N,k)$-words $\ibold$
such that $\pi=\pi_\ibold$ but $\ibold$ is not distinguished does not exceed $C_\pi' N^{k/2}$.
\end{Lemma}
\proof The proof is similar to that of Lemma~\ref{Lemma:CrudeCounting}. We omit the details. \qed

\subsection{Evaluation of $EX^{(N)}_\ibold$ for distinguished $\ibold$}

\subsubsection{Definitions of $\tau_\pi$ and $\sigma_\pi$}
Let $\pi$ be a set partition of $k$. Enumerate $\pi$ and its parts thus:
$$\pi=\{I_1,\dots,I_{|\pi|}\},\;\;\;\min I_1<\dots<\min I_{|\pi|},$$
$$I_\alpha=\{i_{\alpha 1}<\dots<i_{\alpha ,|I_\alpha|}\}\;\mbox{for $\alpha=1,\dots,|\pi|$}.$$
Put
$$\tau_\pi=(i_{11}\cdots i_{1,|I_1|})\cdots
(i_{|\pi|, 1}\cdots i_{|\pi|,|I_{|\pi|}|})\in S_k,\;\;\sigma_\pi=\eta_k^{-1}\tau_\pi\in S_k.$$
By construction
$$\pi(\tau_\pi(i))=\pi(i)$$ for $i=1,\dots,k$.

\begin{Lemma}
\label{Lemma:Involution}
Let $\pi\in \WWW_k$ be a Wigner set partition.
Then the permutation $\sigma_\pi$ is fixed-point-free and squares to the identity. Furthermore, for all distinct $i,j\in \{1,\dots,k\}$, we have
$\{\pi(i),\pi(i+1)\}=\{\pi(j),\pi(j+1)\}\Leftrightarrow
j=\sigma_\pi(i)$.
\end{Lemma}
\proof  Let $G_\pi=(V_\pi,E_\pi)$ be the graph (tree) associated to $\pi$. 
For each $e\in E_\pi$,
there are by Lemma~\ref{Lemma:CrucialInsight}(ii)
exactly two indices $i$ such that
$e=\{\pi(i),\pi(i+1)\}$,
and $\sigma_\pi$ swaps them  by Lemma~\ref{Lemma:CrucialInsight}(iii).
\qed

\begin{Lemma}\label{Lemma:Push}
Let $\pi\in \WWW_k$ be a Wigner set partition.
Put $\sigma=\sigma_\pi$ and $\tau=\tau_\pi$.
Then we have
$$EX^{(N)}_\ibold=
\prod_{\begin{subarray}{c}
\alpha\in \{1,\dots,k\}\\
\mathrm{s.t.}\,\alpha\leq \sigma(\alpha)
\end{subarray}}EX^{(N)}_{\ibold(\alpha),\ibold(\tau(\sigma(\alpha)))}
X^{(N)}_{\ibold(\sigma(\alpha)),\ibold(\tau(\alpha))}
$$
for every $(N,k)$-word $\ibold$ such that $\pi_\ibold=\pi$.
\end{Lemma}
\proof Put $\eta=\eta_k$.  By definition we have $\eta\sigma=\tau$. By the previous lemma we have $\tau\sigma=\eta$. Therefore after rearranging the product $$X_\ibold^{(N)}=\prod_{\alpha\in \{1,\dots,k\}}
X_{\ibold(\alpha),\ibold(\eta(\alpha))},$$ we have
$$
X^{(N)}_\ibold=
\prod_{\begin{subarray}{c}
\alpha\in \{1,\dots,k\}\\
\mathrm{s.t.}\,\alpha\leq \sigma(\alpha)
\end{subarray}}X^{(N)}_{\ibold(\alpha),\ibold(\tau(\sigma(\alpha)))}
X^{(N)}_{\ibold(\sigma(\alpha)),\ibold(\tau(\alpha))}.
$$
It suffices to prove enough independence to justify  pushing the expectation under the product.
For every $\gamma\in \{1,\dots,k\}$ let $i(\gamma)\leq j(\gamma)$
be the integers $\ibold(\gamma)$ and $\ibold(\eta(\gamma))$ rearranged.
By definition of $\pi_\ibold$ and the preceding lemma,
for all distinct 
$\alpha,\beta\in \{1,\dots,k\}$, if 
$\max(|i(\alpha)-i(\beta)|,|j(\alpha)-j(\beta)|)\leq K$,
then $\beta=\sigma(\alpha)$.
By Assumption~\ref{Assumption:Model}(IIb) 
we deduce the desired independence. \qed

\subsubsection{The difference operator 
associated to a set partition $\pi$}
Let $\pi$ be a set partition of $k$. Let $\tau=\tau_\pi\in S_k$
be the canonically associated permutation.
Let $f:\{1,\dots,k\}\rightarrow \ZZ$ be a function.
We define $\partial_\pi f:\{1,\dots,k\}\rightarrow \ZZ$
by the formula
$(\partial_\pi f)(i)=f(i)-f(\tau(i))$
for $i=1,\dots,k$.
Given a function $g:\{1,\dots,k\}\rightarrow \ZZ$,
the equation $\partial_\pi f=g$ 
has a solution $f:\{1,\dots,k\}\rightarrow \ZZ$ if and only if $\sum_{\alpha\in A}g(\alpha)=0$
for all parts $A\in \pi$, and $f$ is unique
up to the addition of a function $\{1,\dots,k\}\rightarrow \ZZ$
constant on every part of $\pi$.

\begin{Lemma}\label{Lemma:KeyEvaluation}
Let $\pi\in \WWW_k$ be a Wigner set partition.
Put $\sigma=\sigma_\pi$, $\tau=\tau_\pi$ 
and $\partial=\partial_\pi$. Fix $\alpha\in \{1,\dots,k\}$
such that $\alpha\leq \sigma(\alpha)$.
We have
\begin{equation}\label{equation:KeyEvaluation}
\textstyle E X^{(N)}_{\ibold(\alpha),\ibold(\tau(
\sigma(\alpha)))}X^{(N)}_{\ibold(\sigma(\alpha)),\ibold(
\tau(\alpha))}
=s_{(\partial\ibold)(\alpha),(\partial\ibold)(\sigma(\alpha))}
\left(\frac{\ibold(\min\pi(\alpha))}{N},\frac{\ibold(
\min(\pi(\sigma(\alpha))))}{N}\right)
\end{equation}
for every distinguished $(N,k)$-word $\ibold$
such that $\pi_\ibold=\pi$.
\end{Lemma}
\noindent In particular, 
the expectation in question vanishes unless $|\partial \ibold|\leq K$ by Assumption~\ref{Assumption:Model}(IIa).
\proof 
 Let $E(\alpha)$ be the left side of (\ref{equation:KeyEvaluation}).
Assume at first that 
$\ibold(\alpha)\leq \ibold(\tau(\sigma(\alpha))$.
Then $\ibold(\alpha)+3Kk<\ibold(\tau(\sigma(\alpha))$
by property (V) of a distinguished $(N,k)$-word,
and hence
$\ibold(\tau(\alpha))+Kk<\ibold(\sigma(\alpha))$
by property (IV) of a distinguished $(N,k)$-word.
So we have
\begin{equation}\label{equation:FirstTry}
\textstyle
E(\alpha)=
EX_{\ibold(\alpha),\ibold(\tau(\sigma(\alpha)))}^{(N)}\overline{X}_{\ibold(\tau(\alpha)),\ibold(\sigma(\alpha))}^{(N)}=s_{(\partial\ibold)(\alpha),
(\partial \ibold)(\sigma(\alpha))}\left(
\frac{\ibold(\alpha)}{N},
\frac{\ibold(\tau(\sigma(\alpha)))}{N}\right)
\end{equation}
by Assumption~\ref{Assumption:Model}(IIc) and
property (VI) of a distinguished $(N,k)$-word.
Then we deduce that the desired formula holds by Assumption~\ref{Assumption:Kernel}(III)
and property (VII) of a distinguished $(N,k)$-word. Assume next and finally that $\ibold(\alpha)\geq \ibold(\tau(\sigma(\alpha))$.
Reasoning as above we have
$$\textstyle E(\alpha)=EX^{(N)}_{\ibold(\sigma(\alpha)),\ibold(\tau(\alpha))}\overline{X}^{(N)}_{\ibold(\tau(\sigma(\alpha)),\ibold(\alpha)}=
s_{
(\partial \ibold)(\sigma(\alpha)),(\partial\ibold)(\alpha)}
\left(\frac{\ibold(\sigma(\alpha))}{N},
\frac{\ibold(\tau(\alpha))}{N}\right).
$$
Now we apply the symmetry (\ref{equation:Symmetries}),
 and then continue to 
reason as above.
We deduce the desired formula just as before. \qed

\subsection{Contribution of $EX^{(N)}_\ibold$
to the limit for distinguished $\ibold$}
\begin{Lemma}
Let $\pi\in \WWW_k$ and an $(N,k)$-word $\ibold$
be given such that the following hold:
\begin{itemize}
\item For all $A\in \pi$ we have $\ibold(\min A)\in Q^{(N)}_{(k+1)K}$.
\item For all distinct $A,B\in \pi$ we have
$|\ibold(\min A)-\ibold(\min B)|>5Kk$.
\item $|\partial \ibold|\leq K$.
\end{itemize}
Then $\pi_\ibold=\pi$ and (hence) $\pi$ is distinguished.
\end{Lemma}
\noindent No proof is needed, but this point bears emphasis
as an important step in the proof of Proposition~\ref{Proposition:MainLimit}.

\begin{Lemma}\label{Lemma:BigKahuna}
Fix a Wigner set partition $\pi\in \WWW_k$. Then we have
\begin{equation}\label{equation:BigKahuna}
\E M_\pi= \lim_{N\rightarrow\infty}
N^{-k/2-1}\sum_{\ibold}EX^{(N)}_{\ibold},
\end{equation}
where the sum is extended over  distinguished $(N,k)$-words $\ibold$
such that $\pi_\ibold=\pi$.
\end{Lemma}
\proof  
Let $\partial=\partial_\pi$,
$\sigma=\sigma_\pi$, and
$\tau=\tau_\pi$. Let
 $\{t_A\}$ (resp., $\{z_A\}$) be an i.i.d.\ family of random variables
uniform in $[0,1]$ (resp., $S^1$), indexed by finite nonempty 
sets $A$ of positive integers.
We further suppose that the families $\{t_A\}$ and 
$\{z_A\}$ are defined on a common probability space and are independent. 
We denote expectations with respect to these variables
by ${\bf E}$.
We have by Lemma~\ref{Lemma:Involution} and the definitions that
$$\begin{array}{rcl}
\E M_\pi&=&\displaystyle
{\bf E}\prod_{
\begin{subarray}{c}
i\in \{1,\dots,k\}\\
\mathrm{s.t.}\,i\leq \sigma(i)
\end{subarray}}\left(
\sum_{m\in \ZZ}\sum_{n\in \ZZ}s_{mn}(t_{\pi(i)},t_{\pi(\sigma(i))})
z_{\pi(i)}^m z_{\pi(\sigma(i))}^n\right)\\\\
&=&\displaystyle
\sum_{f:\{1,\dots,k\}\rightarrow \ZZ}
{\bf E}
\prod_{
\begin{subarray}{c}
i\in \{1,\dots,k\}\\
\mathrm{s.t.}\,i\leq \sigma(i)
\end{subarray}}
s_{f(i),f(\sigma(i))}(t_{\pi(i)},t_{\pi(\sigma(i))})
\cdot \prod_{A\in \pi}
z_A^{\sum_{i\in A}f(i)}\\\\
&=&\displaystyle\sum_{\begin{subarray}{c}
f:\{1,\dots,k\}\rightarrow \ZZ\cap [-K,K]\\
\forall A\in \pi,\;\sum_{i\in A}f(i)=0
\end{subarray}}
{\bf E}\prod_{
\begin{subarray}{c}
i\in \{1,\dots,k\}\\
\mathrm{s.t.}\,i\leq \sigma(i)
\end{subarray}}
s_{f(i),f(\sigma(i))}(t_{\pi(i)},t_{\pi(\sigma(i))}).
\end{array}
$$
At the last step we integrate out the $z$'s and take into 
account Assumption \ref{Assumption:Model}(IIa).  
We then have
$$\E M_\pi
=\sum_{\begin{subarray}{c}
f:\{1,\dots,k\}\rightarrow\ZZ\\
|\partial f|\leq K\\
\forall A\in \pi\;f(\min A)=f_0(A)
\end{subarray}}{\bf E}\prod_{
\begin{subarray}{c}
i\in \{1,\dots,k\}\\
\mathrm{s.t.}\,i\leq \sigma(i)
\end{subarray}}
s_{(\partial f)(i),
(\partial f)(\sigma(i))}(t_{\pi(i)},t_{\pi(\sigma(i))}),
$$
where $f_0:\pi\rightarrow \ZZ$ is any fixed 
function defined on the parts of $\pi$.
After some straightforward bookkeeping which we omit, it follows by the preceding two lemmas that (\ref{equation:BigKahuna})
holds with the summation on $\ibold$
restricted to distinguished $(N,k)$-words
such that $\pi_\ibold=\pi$ and $|\partial \ibold|\leq K$.
But then by Assumption~\ref{Assumption:Model}(IIa), the limit does not change if we drop the restriction $|\partial\ibold|\leq K$ (the further terms all vanish), whence the result.
\qed

\subsection{Completion of the proof of Proposition~\ref{Proposition:MainLimit}}
\label{subsection:MainLimit}
The suite of lemmas proved in \S\ref{subsection:Suite} shows that the limit on the left side of equation (\ref{equation:Promise}) does not change if we restrict attention to distinguished $(N,k)$-words. 
The last lemma above evaluates the limit
on the left side of (\ref{equation:Promise})
with $\ibold$ restricted to distinguished $(N,k)$-words, and gives the desired value.
The proof of Proposition~\ref{Proposition:MainLimit} is complete.
 \qed

\section{Completion of the proof of Theorem~\ref{Theorem:MainResult}}
\label{section:BigCompletion}
Fix a positive integer $k$. As in \cite[pf.\ of Thm.~3.2, Section 6, p.\ 305]{AZ}, 
Theorem \ref{Theorem:MainResult} will follow as
soon as we can prove that 
\begin{equation}
\label{eq-190606a}
 \mbox{\rm Var}( \langle L^{(N)}, x^k\rangle)\to_{N\to\infty } 0\,.
\end{equation}
We will prove this
by lightly modifying the arguments of \S\ref{subsection:Suite}. 
Now 
\begin{equation}
\label{eq-190606b}
 \mbox{\rm Var}( \langle L^{(N)}, x^k\rangle)=
N^{-k-2} \sum_{(\ibold,\jbold)}\left(E[X_{\ibold}^{(N)}
X_{\jbold}^{(N)}]-E[X_{\ibold}^{(N)}]E[X_{\jbold}^{(N)}]
\right)\,
\end{equation}
where the sum is extended over pairs $(\ibold,\jbold)$ of $(N,k)$-words. 
By Assumption~\ref{Assumption:Model}(Ib) and the H\"{o}lder inequality,
it is (more than) enough to show that the number of pairs
$(\ibold,\jbold)$ making a nonzero contribution to the sum on the right side of
(\ref{eq-190606b}) is $O(N^{k+1})$.
Now fix a pair $(\ibold,\jbold)$ of $(N,k)$-words
indexing a nonzero term in the sum on the right side of
(\ref{eq-190606b}).
Let $\ibold\jbold$ be the $(N,2k)$-word obtained
by concatenating $\ibold$ and $\jbold$,
i.~e., 
$$\ibold\jbold(\alpha)=
\left\{\begin{array}{rl}
\ibold(\alpha)&\mbox{if $\alpha\in \{1,\dots,k\}$,}\\
\jbold(\alpha-k)&\mbox{if $\alpha\in \{k+1,\dots,2k\}$.}
\end{array}\right.
$$
Consider the set partition $\pi=\pi_{\ibold\jbold}$ defined in \S\ref{subsubsection:BlobDef}. We need also to consider a graph
associated to $\pi$. The relevant graph
is no longer $G_\pi$, but rather a slightly modified version $\tilde{G}_\pi=(\tilde{V}_\pi,\tilde{E}_\pi)$, where
$\tilde{V}_\pi=\pi$ and
$$\tilde{E}_\pi=\{\{\pi(1),\pi(2)\},\dots,\{\pi(k),\pi(1)\}\}
\cup\{\{\pi(k+1),\pi(k+2)\},\dots,\{\pi(2k),\pi(k+1)\}\}.$$
By construction $\tilde{G}_\pi$ comes equipped with
two walks, namely
$$\pi(1),\dots,\pi(k),\pi(1)\;\;\mbox{and}\;\;\pi(k+1),\dots,\pi(2k),\pi(k+1).$$
Arguing as in the proof of Lemma~\ref{Lemma:WignerBound},
we find that nonvanishing of the term on the right side of (\ref{eq-190606b})
indexed by $(\ibold,\jbold)$ implies that the walks
jointly visit every edge of $\tilde{G}_\pi$ at least twice,
and moreover there must exist some edge of $\tilde{G}_\pi$
visited by both walks. Thus $|\tilde{E}_\pi|\leq k$ and
$\tilde{G}_\pi$ is connected.
It follows that $|\pi|=|\tilde{V}_\pi|\leq k+1$.
Finally, arguing as in the proof of Lemma~\ref{Lemma:CrudeCounting},
we find that the number of nonzero terms 
on the right side of (\ref{eq-190606b}) is indeed $O(N^{k+1})$.
The proof of Theorem~\ref{Theorem:MainResult} is complete. \qed

\section{An algebraicity criterion}\label{section:Soft}
In this section, which is completely independent of the preceding ones, we develop a ``soft'' method for proving that a holomorphic function is algebraic under hypotheses commonly encountered in random matrix theory.

\subsection{Formulation of an algebraicity criterion}
\subsubsection{Notation}
Let $\{X_i\}_{i=1}^\infty$ be independent (algebraic) variables.
Let \linebreak $\CC[X_1,\dots,X_n]$
denote the ring of polynomials in $X_1,\dots,X_n$
with coefficients in $\CC$. We view $\CC[X_1,\dots,X_n]$
as a subring of $\CC[X_1,\dots,X_{n+1}]$.
Given $F\in \CC[X_1,\dots,X_n]$, let $F(0)\in \CC$
be the result of setting $X_1=\cdots=X_n=0$.
Let $\CC(X_1,\dots,X_n)$ be the field of rational functions
in the variables $X_1,\dots,X_n$, i.~e., the ring consisting
of fractions $A/B$ where $A,B\in \CC[X_1,\dots,X_n]$ and $B$ does not vanish identically.
We say that $F\in \CC(X_1,\dots,X_n)$ is {\em defined at the origin}
if $F=A/B$ with $A,B\in \CC[X_1,\dots,X_n]$ such that $B(0)\neq 0$,
in which case we put $F(0)=A(0)/B(0)$.

\subsubsection{DIRE families}
Let $\varphi_1,\dots,\varphi_N$ be a finite family of
holomorphic functions each defined in a 
connected open neighborhood of the origin in $\CC^n$;
the domains need not be the same.
We will say that $\varphi_1,\dots,\varphi_N$ are {\em defined implicitly 
by rational equations} ({\em DIRE}  for short)
if there exist $\Phi_1,\dots,\Phi_N\in \CC(X_1,\dots,X_{n+N})$ such that
\begin{enumerate}
\item $\Phi_i$ is defined at the origin 
and $\Phi_i(0)=0$ for $i=1,\dots,N$,
\item $(\det_{i,j=1}^N \frac{\partial 
\Phi_i}{\partial X_{j+n}})(0)\neq 0$, and
\item $\Phi_i(z_1,\dots,z_n,\varphi_1(z)-\varphi_1(0),\dots,
\varphi_N(z)-\varphi_N(0))=0$
for $i=1,\dots,N$ and \linebreak 
$z=(z_1,\dots,z_n)\in \CC^n$ sufficiently near the origin.
\end{enumerate}
The relationship between this definition and 
the implicit function theorem
for holomorphic functions \cite[Proposition 6.1]{cartan}
is close. Indeed,
given  $\Phi_1,\dots,\Phi_N\in \CC(X_1,\dots,X_{n+N})$
with properties (I,II) above, the implicit function theorem for holomorphic functions says that there exist
holomorphic functions $\varphi_1,\dots,\varphi_N$ each defined in a connected open neighborhood of the origin such that (III) holds,
and the theorem further asserts uniqueness of these functions
in the sense that if $\psi_1,\dots,\psi_N$ are holomorphic functions each defined in a connected open neighborhood of the origin in $\CC^n$ and also satisfying (III), then
for $i=1,\dots,N$ there exists a neighborhood of the origin
in which $\psi_i$ and $\varphi_i$ 
differ by a constant.

In the sequel, for brevity, 
when we say ``$\varphi_1,\dots,\varphi_N$ is an $n$-variable
DIRE family'', this is short for the assertion that
``$\varphi_1,\dots,\varphi_N$ are holomorphic functions each defined in some connected open neighborhood of the origin in $\CC^n$ which together have the property of being defined implicitly by rational equations''.

\subsubsection{Algebraic functions}
A holomorphic function $\varphi$ defined in a nonempty open subset $D\subset \CC^n$ is called an {\em $n$-variable algebraic function} if there exists a not-identically-vanishing polynomial $F=F(X_1,\dots,X_{n+1})\in \CC[X_1,\dots,X_{n+1}]$
such that $F(z_1,\dots,z_n,\varphi(z))=0$
for all $z=(z_1,\dots,z_n)\in D$. 
If $D$ is connected
and $U\subset D$ is any nonempty open subset, then algebraicity
of $F$ in $U$ implies algebraicity 
of $F$ in $D$ by the principle of analytic continuation.\\

The main result of Section \ref{section:Soft} is the following.

\begin{Theorem}\label{Theorem:BasicFunctionalEquation}
Let $\varphi_1,\dots,\varphi_N$ be an $n$-variable DIRE family.
Then each $\varphi_i$ is an $n$-variable algebraic function.
\end{Theorem}

\noindent The proof takes up the last
 several subsections of Section \ref{section:Soft}. 
Before turning to the proof we give key examples of DIRE families, and describe techniques for constructing new DIRE families from old.

\begin{Proposition}\label{Proposition:RandomWalk}
Fix a positive integer $\ell$ and put $L=2\ell+1$.
For $z=(z_1,\dots,z_L)\in \CC^L$
such that $\sum_{i=1}^L |z_i|<1$
and integers $j=1,\dots,L$,
consider the quantity
\begin{equation}\label{equation:Thetas}
\vartheta_j(z)=-\delta_{j,\ell+1}+\frac{1}{2\pi}\int_0^{2\pi}
\frac{\exp(-{\bf i} (j-\ell-1)x)dx}{1-\sum_{k=1}^L
z_k\exp({\bf i} (k-\ell-1)x)},
\end{equation}
where ${\bf i}^2=-1$,
which depends holomorphically on $z$
and vanishes for $z=0$.
Then the family $\vartheta_1,\dots,\vartheta_L$
can be extended to an $L$-variable DIRE family
$\vartheta_1,\dots,\vartheta_N$.\end{Proposition}
\proof
Fix $z=(z_1,\dots,z_{L})\in \CC^{L}$ such that $\sum_{i=1}^L|z_i|<1$.
Let  $p=[p_{ij}]_{i,j\in \ZZ}$
be the unique matrix of complex numbers with rows and columns indexed by $\ZZ$ with the following properties:
\begin{itemize}
\item $p_{ij}=0$ for all $i$ and $j$ such that  $|i-j|>\ell$.
\item $p_{ij}=z_{j-i+\ell+1}$ for all $i$ and $j$ such that
$|i-j|\leq \ell$.
\end{itemize}
Let $|p|$ be the $\ZZ$-by-$\ZZ$ matrix with entries $|p_{ij}|$.
We have the crude estimate
\begin{equation}\label{equation:pMajorization}
(|p|^t)_{ij}\leq \left(\sum_{i=1}^L |z_i|\right)^t
\end{equation}
holding for all positive integers $t$.

Now supposing $p$ were a Markov matrix (which of course it is not),
we could view each entry $p_{ij}$ as a transition probability
$$P(S_{t+1}=j\mid S_t=i)=p_{ij},$$
for a random walk
$\{S_t\}_{t=0}^\infty$ on $\ZZ$ with step-length bounded by $\ell$ and we would have, for any subset $M\subset \ZZ^{t}$,
an equality
$$P((S_1,\dots,S_t)\in M\mid S_0=i_0)=
\sum_{(i_1,\dots,i_t)\in M}p_{i_0i_1}\cdots p_{i_{t-1}i_t}.$$
Because it is a valuable aid to intuition, we will
make the line above a definition.  Notice that by (\ref{equation:pMajorization}) the sum on the right is absolutely convergent, and that the sum of the absolute values of the terms is $\leq (\sum_{i=1}^{L}|z_i|)^t$.
We will be able to calculate using the usual rules of probability, provided we never invoke
 positivity $p_{ij}\geq 0$ or the Markov property
$\sum_j p_{ij}=1$. 
We are interested in functions of $z$
describable in the language of random walk because, as we explain presently, the functions $\vartheta_j$ belong to this class, and moreover within this class we can easily find the ``extra'' functions needed to extend $\vartheta_1,\dots,\vartheta_{L}$ to a DIRE family.

Let $U,V,A,B,C$ be $\ell$-by-$\ell$ matrices with complex entries defined as follows:
$$\begin{array}{rcl}U_{ij}&=&
\sum_{t=1}^\infty P(S_t=j,\min_{u=0}^t S_u>0\vert S_0=i)
,\\
V_{ij}&=&\sum_{t=1}^\infty P(S_t=j,\max_{u=0}^t S_u<\ell+1\vert S_0=i),\\
A_{ij}&=&P(S_1=j+\ell\mid S_0=i),\\
B_{ij}&=&P(S_1=j\mid S_0=i),\\
C_{ij}&=&P(S_1=j-\ell\mid S_0=i).
\end{array}
$$
Let $W,D$ be $L$-by-$L$ matrices
with complex entries defined as follows:
$$
\begin{array}{rcl}
W_{ij}&=&\sum_{t=1}^\infty P(S_t=j\vert S_0=i),\\
D_{ij}&=&P(S_1=j\mid S_0=i).
\end{array}
$$
By (\ref{equation:pMajorization}), all the series
in question here converge absolutely
and define functions of $z$ holomorphic in the domain
$\sum_{i=1}^L|z_i|<1$. Notice that $A,B,C,D$ depend linearly on $z$,
and that all matrices $A,B,C,D,U,V,W$
vanish for $z=0$.

Consider the expansion of the integrand of (\ref{equation:Thetas}) by geometric series and then integrate term by term. One finds in this way a series expression for $\vartheta_j$ identical to the series expression
defining some entry of the matrix $W$. In short, every $\vartheta_j$
appears in $W$.

By breaking paths down according to visits to the sets $\{-\ell,\dots,-1\}$, $\{0\}$ and $\{1,\dots,\ell\}$,
we obtain in the usual way recursions
\begin{equation}
\label{equation:TheRecursions}
\begin{array}{rcl}
0&=&-U+B+A(1+U)C(1+U),\\
0&=&-V+B+C(1+V)A(1+V),\\
0&=&-W+D+\left[\begin{array}{ccc}
C(1+V)A&0&0\\
0&0&0\\
0&0&A(1+U)C
\end{array}\right](1+W).
\end{array}
\end{equation}
Extend the given family $\vartheta_1,\dots,\vartheta_n$
to an enumeration $\vartheta_1,\dots,\vartheta_N$
of all entries of $U$, $V$ and $W$. Rewrite the system of equations (\ref{equation:TheRecursions}) as a system of $N$ polynomial equations
in $z_1,\dots,z_{L},\vartheta_1,\dots,\vartheta_N$,
in order to find polynomials $$\Phi_i(X_1,\dots,X_{L+N})\in \CC[X_1,\dots,X_{L+N}]\;\;\mbox{for}\; i=1,\dots,N$$
such that 
$$\Phi_i(0,\dots,0,X_{L+1},\dots,X_{L+N})=-X_{i+L},\;\;
\Phi_i(z_1,\dots,z_L,\vartheta_1,\dots,\vartheta_N)=0.$$
Note that part (II) of the definition of DIRE family
holds because $\frac{\partial \Phi_i}{\partial X_{j+L}}(0)=-\delta_{ij}$.
Thus $\vartheta_1,\dots,\vartheta_N$
is indeed a DIRE family.
\qed

\subsection{Natural operations on DIRE families}
\label{subsection:NaturalOperations}
We write down some lemmas which will be helpful in applying
the notion of DIRE family. The first three are trivial
but deserve being stated for the sake of emphasis.
The last is the trick decisive for the application of Theorem~\ref{Theorem:BasicFunctionalEquation} to the proof of Theorem~\ref{Theorem:MainResultBis}.

\begin{Lemma}\label{Lemma:Op1}
Let $\varphi_1,\dots,\varphi_N$ and $\psi_1,\dots,\psi_M$
be $n$-variable DIRE families.
Then the concatenation $\varphi_1,\dots,\varphi_N,\psi_1,\dots,\psi_M$
is  an $n$-variable DIRE family.
\end{Lemma}
\begin{Lemma}\label{Lemma:Op2}
Let $\varphi_1,\dots,\varphi_N$ be an $n$-variable DIRE family.
Let $\varphi_{N+1}$ be a $\CC$-linear combination of $\varphi_1,\dots,\varphi_N$. Then the extended family $\varphi_1,\dots,\varphi_{N+1}$
is an $n$-variable DIRE family.
\end{Lemma}
\begin{Lemma}\label{Lemma:Op3}
Let $\varphi_1,\dots,\varphi_N$ be an $n$-variable DIRE family.
Let $A$ be an $n$ by $m$ matrix with complex entries.
Identify $\CC^n$ and $\CC^m$ with spaces of column vectors.
Then there exists an $m$-variable DIRE family $\psi_1,\dots,\psi_N$
such that for $i=1,\dots,N$ we have $\psi_i(z)=\varphi_i(Az)$
for all $z\in \CC^m$ sufficiently near the origin.
\end{Lemma}

\begin{Lemma}\label{Lemma:McGuffin}
Let $\psi_1,\dots,\psi_n$ be a family of holomorphic
functions defined in an open disk centered at the origin in $\CC$.
Let $\varphi_1,\dots,\varphi_N$ be an $n_0$-variable DIRE family,
where $N\geq n\geq n_0$.
Assume that
$$\psi_i(z)=z\varphi_i(z\psi_1(z),\dots,z\psi_n(z))$$
for $i=1,\dots,n$ and $z\in \CC$ sufficiently near the origin. Then $\psi_1,\dots,\psi_n$ can be extended to a $1$-variable DIRE family.
\end{Lemma}
\proof By the preceding lemma we may assume without loss of generality that $n=n_0$.
For $i=n+1,\dots,n+N$ the formula
$$\psi_{i}(z)=\varphi_{i-n}(z\psi_1(z),\dots,
z\psi_n(z))-\varphi_{i-n}(0)$$
defines a holomorphic function $\psi_{i}$ in some open neighborhood of the origin in $\CC$.
We will prove that the extended family $\psi_1,\dots,\psi_{N+n}$ is a $1$-variable DIRE family.
Note that all the functions $\psi_i$ vanish at the origin.
 Let $\Phi_1,\dots,\Phi_{n+N}\in \CC(X_1,\dots,X_{N+n})$ be with respect to $\varphi_1,\dots,\varphi_N$ as called for by the definition of an $n$-variable DIRE family. For $i=1,\dots,N+n$  
 define $\Psi_i\in \CC(X_1,\dots,X_{N+n+1})$ by the formula
$$
 \begin{array}{cl}
 &\Psi_i(X_1,\dots,X_{n+N+1})\\\\
 =&
 \left\{\begin{array}{lr}
X_{i+1}-X_1(X_{i+n+1}+\varphi_i(0))&\mbox{if $1\leq i\leq n$,}\\
\Phi_{i-n}(X_1X_2,\dots,X_1X_{n+1},X_{n+2},\dots,X_{n+N+1})&
\mbox{if $n+1\leq i\leq n+N$.}
 \end{array}\right.
 \end{array}
 $$
Then:
 \begin{enumerate}
\item $\Psi_i$ is defined at the origin 
and $\Psi_i(0)=0$ for $i=1,\dots,N+n$,
\item $(\det_{i,j=1}^{N+n} \frac{\partial 
\Psi_i}{\partial X_{j+1}})(0)\neq 0$, and
\item $\Psi_i(z,\psi_1(z),\dots,
\psi_{N+n}(z))=0$
for $i=1,\dots,N+n$ and 
$z\in \CC$ sufficiently near the origin.
\end{enumerate}
In other words,
$\Psi_1,\dots,\Psi_{n+N}$
are with respect to 
$\psi_1,\dots,\psi_{n+N}$ as called for by the definition
of $1$-variable DIRE family. 
\qed

 We next formulate a purely algebraic result
and explain how to deduce 
Theorem~\ref{Theorem:BasicFunctionalEquation}
from it.

\begin{Theorem}\label{Theorem:Algebraicity}
Let $n$ and $N$ be positive integers.
Let $$F_1,\dots,F_{n+N}\in \CC[X_1,\dots,X_{N+n}]$$
be given with the following properties:
\begin{equation}\label{equation:PreJacCon}
F_i(0)=0\;\;\mbox{for $i=1,\dots,N$}.
\end{equation}
\begin{equation}\label{equation:JacCon}
\left(\det_{i,j=1}^N \frac{\partial F_i}{\partial X_{j+n}}\right)(0)\neq 0.
\end{equation}
Let $I\subset \CC[X_1,\dots,X_{N+n}]$ be the ideal generated by $F_1,\dots,F_N$.
Then there exist polynomials
$G\in \CC[X_1,\dots,X_{n+1}]$ and 
$H\in \CC[X_1,\dots,X_{n+N}]$ such that
$G\neq 0$, $H(0)\neq 0$ and $GH\in I$.
\end{Theorem}
\noindent \textbf{Remark}. The important point here is that
$G$ is a polynomial involving only the variables $X_1,\dots,X_{n+1}$;
the variables $X_{n+2},\dots,X_{n+N}$ are uninvolved. 
The proof of  Theorem~\ref{Theorem:Algebraicity} is a routine
application of the
theory of noetherian local rings, and will be given
in \S
\ref{subsection:AlgebraicityProof} after
we review in \S \ref{subsection:Review}
the needed material from commutative algebra.
\subsection{Deduction of 
Theorem~\ref{Theorem:BasicFunctionalEquation}  from
Theorem~\ref{Theorem:Algebraicity}}
By symmetry it
 is enough to show that $\varphi_1$ is algebraic.
Let $\Phi_1,\dots,\Phi_N\in \CC(X_1,\dots,X_{n+N})$
be as required to exist by the definition of $n$-variable DIRE family with respect to $\varphi_1,\dots,\varphi_N$.
Write 
$$\Phi_i=F_i/D_i\;\;\;(F_i,D_i\in \CC(X_1,\dots,X_{n+N}),\;D_i(0)\neq 0).$$
Without loss of generality we may simply assume that $D_i=1$ and hence $\Phi_i=F_i$.
Then conditions (I,II) are precisely the hypotheses (\ref{equation:PreJacCon},\ref{equation:JacCon})
of Theorem~\ref{Theorem:Algebraicity}.
Let $G$ and $H$ be as provided by Theorem~\ref{Theorem:Algebraicity}.
Then $H(z_1,\dots,z_n,\varphi_1(z),\dots,\varphi_N(z))$
is nonvanishing for $z=(z_1,\dots,z_n)\in \CC^n$
sufficiently near the origin,
hence \linebreak $G(z_1,\dots,z_n,\varphi_1(z))$ vanishes 
for $z\in \CC^n$ sufficiently near the origin, and hence
$\varphi_1$ is indeed algebraic.
\qed

\subsection{Review of dimension theory of noetherian local rings}
\label{subsection:Review}
The material reviewed here is developed in detail 
in the texts \cite{AtiyahMacdonald} and \cite{Matsumura}. 


\subsubsection{The setting}\label{subsubsection:Flag0}
In our review {\em rings} are always commutative with unit.
A ring $R$ is {\em noetherian} if every ideal is finitely generated. 
All rings to be considered below are assumed to be noetherian.
A {\em local ring} is a ring possessing
a unique maximal ideal. 
For the rest of the discussion we fix a noetherian local ring $R$ with  maximal ideal $M$ and denote
the residue field $R/M$ by $k$. We urge 
the reader to keep the following key example of triples 
$(R,M,k)$ in mind:
\begin{itemize}
\item $R=\{F\in \CC(X_1,\dots,X_d)\mid 
\mbox{$F$ is defined at the origin}\},$
\item $M=\{F\in R\mid F(0)=0\}$, and
\item $k=\CC$.
\end{itemize}

\subsubsection{Dimension of a noetherian local ring}
\label{subsubsection:Flag1}
For each integer $n$ the quotient \linebreak $M^n/M^{n+1}$
is a finite-dimensional vector space over $k$. (Here $M^n$ stands for the ideal generated by all $n$-fold products of elements of $M$, and by convention $M^0=R$.)
Consider the nonnegative-integer-valued function
$$\chi(n)=\sum_{i=0}^{n-1}\dim_k M^i/M^{i+1}$$
of nonnegative integers $n$.
There exists a unique polynomial $F(t)$ in a variable $t$ with rational coefficients such that $\chi(n)=F(n)$ for all $n\gg 0$. 
(See \cite[12.C]{Matsumura} or \cite[Cor.\ 11.5]{AtiyahMacdonald}.)
The polynomial $F(t)$ is called the
 {\em Hilbert-Samuel} polynomial of $R$. The degree in $t$
of $F(t)$ is by definition the {\em dimension} of $R$,
and denoted $d(R)$. In the example discussed in 
\S \ref{subsubsection:Flag0}, $d(R)=d$.
\subsubsection{Regular local rings}\label{subsubsection:Flag2}
We say that $R$ is {\em regular} if 
$d(R)=\dim_k M/M^2$, in which case necessarily $\chi(n)=
\left(\begin{array}{c}n+d(R)\\
d(R)\end{array}\right)$ for all $n$, and $R$ is an integral domain. 
(See \cite[(17.E) Thm.\ 35 and (17.F) Thm.\ 36]{Matsumura} or
\cite[Thm.\ 11.22 and Lemma 11.23]{AtiyahMacdonald}. In the case of
the example in \S \ref{subsubsection:Flag0}, this can be verified
directly, noting that $\dim M/M^2=d$ in that case).
Suppose for the rest of this paragraph that $R$ 
is regular of dimension $d$.
Elements $f_1,\dots,f_d\in M$ forming a basis
over $k$ for the quotient $M/M^2$ 
are said to form a {\em regular system of parameters}
for $R$. By a standard argument employing Nakayama's lemma 
(for the latter see \cite[(1.M) Lemma]{Matsumura} 
or \cite[Prop.\ 2.6]{AtiyahMacdonald}) any 
regular system of parameters for $R$ necessarily
generates the maximal ideal $M$. In the key example of
\S \ref{subsubsection:Flag0},
the variables $X_1,\dots,X_d$ form a regular system of parameters.

\subsubsection{Cutting down regular 
local rings}\label{subsubsection:Flag3}
Again suppose that $R$ is regular of dimension $d$.
Given a regular
system of parameters $f_1,\dots,f_d$ in $R$,
and also given $i=0,\dots,d$, the ideal $(f_1,\dots,f_i)\subset R$
generated by $f_1,\dots,f_i$
is prime and the quotient $R/(f_1,\dots,f_i)$
is a regular local ring in which
the images of $f_{i+1},\dots,f_d$ form a regular
system of  parameters. (See \cite[(17.F) Thm.\ 36]{Matsumura}.)
One should think of this fact as an 
algebraic version of the implicit function theorem.

\subsubsection{Algebraic independence
 of regular systems of parameters}
\label{subsubsection:Flag3a}
Again assume that $R$ is regular of dimension $d$, and further assume 
that $R$ contains a field $k_0$. 
Then every regular system of parameters $f_1,\dots,f_d$ in $R$
is algebraically independent over $k_0$, i.~e., for every polynomial
$F(X_1,\dots,X_d)$ in independent 
variables $X_1,\dots,X_d$ with coefficients in $k_0$,
if $F(f_1,\dots,f_d)=0$, then $F=0$.
(See \cite[Cor.\ 11.21]{AtiyahMacdonald}
or \cite[(20.D) App.\ 1]{Matsumura}.)
In the example of \S \ref{subsubsection:Flag0}, 
we may take $k_0=\CC$.

\subsubsection{Relation of dimension to transcendence degree}
\label{subsubsection:Flag4}
Assume now that $(R,M,k)$ 
is of the form of the example from \S \ref{subsubsection:Flag0}.
Let $P$ be any prime ideal of $R$. The quotient $R/P$ is again a 
noetherian local ring (but maybe not regular). (The ring $R/P$ admits
 interpretation as the local ring at a point, perhaps singular, of an 
algebraic variety in $\CC^d$.)
Let $e$ be the {\em transcendence degree} of $R/P$ over $\CC$,
i.~e., the supremum of the set of integers $m\geq 0$
such that there exist $m$ elements of $R/P$ algebraically independent over $\CC$.
Then we have $e=d(R/P)$.  (See \cite[Thm.\ 11.25]{AtiyahMacdonald}.) One has this equality
whether or not $R/P$ is regular.

\subsection{Proof of Theorem~\ref{Theorem:Algebraicity}}
\label{subsection:AlgebraicityProof}
We are ready to move rapidly through the proof.
We will flag the relevant paragraphs above at each step.
Consider the ring $R\subset \CC(X_1,\dots,X_{n+N})$
consisting of all fractions $A/B$ where $A,B\in \CC[X_1,\dots,X_{n+N}]$
and $B(0)\neq 0$.
Then $R$, see 
\S \ref{subsubsection:Flag2},
 is a regular local ring of dimension $n+N$.
Hypotheses (\ref{equation:PreJacCon},\ref{equation:JacCon}) 
imply that  $X_1,\dots,X_n,F_1,\dots F_N$ 
form a regular system of parameters for $R$.
Let $P$ be the prime ideal of 
$R$ generated by $F_1,\dots,F_N$ and let $x_1,\dots,x_{n+N}$
be the images in the quotient $R/P$ of $X_1,\dots,X_{n+N}$,
respectively.
By \S \ref{subsubsection:Flag3},
the quotient $R/P$ is a regular local ring of dimension $n$
for which $x_1,\dots,x_n$ forms a regular system of parameters.
Necessarily, by \S \ref{subsubsection:Flag3a},
$x_1,\dots,x_n$ are algebraically independent
over $\CC$ and furthermore, see \S \ref{subsubsection:Flag4},
no set of elements of $R/P$ algebraically 
independent over $\CC$ can have cardinality exceeding $n$ 
(we emphasize that equality to $0$ is taken here in $R/P$, not $R$). 
So there exists
$0\neq G=G(X_1,\dots,X_{n+1})\in \CC[X_1,\dots,X_{n+1}]$
such that $G(x_1,\dots,x_{n+1})=0$ and
hence (equivalently) $G\in \CC[X_1,\dots,X_{n+1}]\cap P$. Now every element of $P$ can be written $A/B$ where $A\in I$
and $B\in \CC[X_1,\dots,X_{n+N}]$ is such that $B(0)\neq 0$.
In particular we may write $G=A/B$ in such fashion.
Taking $H=B$, we have $GH=A\in I$, as desired. \qed

\section{Proof of Theorem~\ref{Theorem:MainResultBis}}\label{section:SoftApplied}

Let $\ell$ be a large positive integer.
Let $\{Z_m\}_{m=1}^M$ be an enumeration of all complex-valued functions $f$ on color space of the form
$$f(c)=\one_{I}(x)\xi^j/\sqrt{\mbox{length of $I$}};\;\;(c=(x,\xi)\in C=[0,1]\times S^1,\;\;I\in \III,\;\;j=-\ell,\dots,\ell).$$
Note that $\{Z_m\}_{m=1}^M$ is an orthonormal system in $L^2(C)$.
We suppose $\ell$ is chosen large enough so that we have an expansion
$$s(c,c')=\sum_{i=1}^M\sum_{j=1}^M \rho_{ij}Z_i(c)\overline{Z}_j(c')$$
for some complex constants $\rho_{ij}$. Also write
$$1=\sum_{m=1}^M \rho_mZ_m\;\;\;\left(\rho_m=\int \overline{Z}_m(c)P(dc)\right).$$
With $\Psi(c,\lambda)$ as defined in the color equations (\ref{equation:Color}), put
$$w_m(z)=\int \Psi(c,1/z)\overline{Z}_m(c)P(dc)\;\;\;(m=1,\dots,M),\;\;\;\;w_{M+1}(z)=S(1/z).$$
The functions $w_m(z)$ are defined and holomorphic 
for $|z|$ small and positive, and moreover are $O(|z|)$,
and hence extend to holomorphic functions in a small
disk about the origin which vanish at the origin.
Consider the functions
$$\begin{array}{rcl}
F_m(z)&=&\displaystyle\sum_{j=1}^M\rho_{mj}\int
\frac{\overline{Z}_j(c)P(dc)}
{1-\sum_{m=1}^Mz_mZ_m(c)}\;\;\;(m=1,\dots,M),\\\\
F_{M+1}(z)&=&\displaystyle\sum_{j=1}^M\rho_{j}\int\frac{\overline{Z}_j(c)P(dc)}
{1-\sum_{m=1}^Mz_mZ_m(c)}  
\end{array}
$$
defined
and holomorphic for $z=(z_1,\dots,z_M)\in \CC^M$ sufficiently near the origin.
The family $F_1,\dots,F_{M+1}$
may be extended to a $M$-variable DIRE family
by Proposition~\ref{Proposition:RandomWalk}
along with
Lemmas~\ref{Lemma:Op1}, \ref{Lemma:Op2}, and \ref{Lemma:Op3}.
Now the general color equations can be rewritten in the form
\begin{equation}\label{equation:MasterEquationsTer}
w_m(z)=z
F_m\left(z w_1(z),\dots,
zw_M(z)\right)\;\;\;\mbox{for $m=1,\dots,K+1$ and $|z|$ small.}
\end{equation}
By Lemma~\ref{Lemma:McGuffin} the family
$w_1,\dots,w_M$ may be extended to a $1$-variable DIRE family, and hence each function $w_m$
is algebraic by Theorem~\ref{Theorem:BasicFunctionalEquation}.
Finally, $S(\lambda)=w_{M+1}(1/\lambda)$ is algebraic.
\qed

\end{document}